\documentclass[10pt]{article}

\usepackage{amsmath,amsxtra,amssymb,latexsym, amscd,amsthm}
\usepackage{here}
\usepackage{graphpap}
\usepackage{graphicx}
\usepackage{color}
\usepackage{array}
\usepackage{booktabs}
\usepackage[ruled]{algorithm2e}
\usepackage{algorithmic}
\usepackage{ marvosym }

\newtheorem{theorem}{Theorem}[section]

\newtheorem{lemma}[theorem]{Lemma}

\theoremstyle{definition}

\newtheorem{remark}[theorem]{Remark}

\newtheorem{example}{Example}[section]

\numberwithin{equation}{section}

\def\divv{\text{div}}
\def\<{\langle}
\def\>{\rangle}

\def\N{\mathbb{N}}
\def\limn{\lim_{n\to\infty}}
\def\limsupn{\limsup_{n\to\infty}}
\def\liminfn{\liminf_{n\to\infty}}
\definecolor{purple}{rgb}{0.4, 0.0, 0.4}

\def\en{\enskip}

\begin{document}


\begin{center}
\bf  Determining two coefficients in diffuse optical tomography \\with incomplete and noisy Cauchy data
\end{center}


\centerline{ Tran Nhan Tam Quyen \let\thefootnote\relax\footnote{Institute for Numerical and Applied Mathematics, University of Goettingen,
Lotzestr. 16-18, 37083 Goettingen, Germany (quyen.tran@uni-goettingen.de)}}

\begin{quote}
{\small {\bf Abstract}
In this paper we investigate the non-linear and ill-posed inverse problem of simultaneously  identifying the conductivity and the reaction in diffuse optical tomography with noisy measurement data available  on an accessible part of the boundary. We propose an energy functional method and the total variational regularization combining with the quadratic stabilizing term to formulate the identification problem to a PDEs constrained optimization problem. We show the stability of the proposed regularization method and the convergence of the finite element regularized solutions to the identification in the $L^s$-norm for all $s\in [0,\infty)$ and in the sense of the Bregman distance with respect to the total variation semi-norm.  
To illustrate the theoretical results, a numerical case study is presented
which supports our analytical findings.
}


\end{quote}

\section{Introduction}

Electrical impedance tomography is a noninvasive type of medical imaging, where the tomographic image of the  electrical conductivity, permittivity, and impedance of a body part is desired to infer from  surface electrode measurements. This problem attracted a great deal of attention from many applied scientists in the last decades. For surveys on the subject, we refer the reader to, e.g., \cite{borcea,br03,cheney,dijkstra,
Dobson96,mueller,uhlmann}  and the references given there. 

Mathematically, assume that the electric potential or voltage $u$ in the body $\Omega$ is governed by the equation
\begin{align*}
\nabla \cdot \big(q \nabla u \big)  = 0 \quad\mbox{in}\quad \Omega \subset \mathbb{R}^d, \quad d\ge 2
\end{align*}
with a free source. Here $q = q(x), \enskip x\in \Omega$ is the electrical conductivity which must be identified from some measurements of the state $u$ on the boundary $\partial\Omega$ of the body $\Omega$. In an ideal situation we know all the voltages $u_{|\partial\Omega} := g \in H^{1/2}(\partial\Omega)$ and the outward pointing normal component of the current densities $\Lambda^{\partial\Omega}_q g := {q \frac{\partial u}{\partial\vec{n}}}_{|\partial\Omega} := j \in H^{-1/2}(\partial\Omega)$ as well, i.e. the knowledge of the Dirichlet-to-Neumann map
$$\Lambda^{\partial\Omega}_q: H^{1/2}(\partial\Omega)
\rightarrow H^{-1/2}(\partial\Omega)$$
is described, where $\vec{n}$ is the unit outward normal on $\partial\Omega$. This is the continuum model which is commonly used  in mathematical researches on the question of the solution uniqueness 
\begin{align*}
p, q \in \mathcal{C} \subset L^\infty(\Omega) \quad\mbox{with}\quad \Lambda^{\partial\Omega}_p = \Lambda^{\partial\Omega}_q
\quad \Rightarrow \quad p=q.
\end{align*}
In dimensions three and higher the uniqueness result has been investigated by Sylvester and Uhlmann \cite{sylvester}, P\"aiv\"arinta {\it el al.} \cite{paivarinta}, and Brown and Torres \cite{brown}, depending on the smoothness of considered conductivities. Meanwhile for the two dimensional setting it can be found in Nachman \cite{nachman}, Brown and Uhlmann \cite{brown2}, and Astala and P\"aiv\"arinta \cite{aspa06}.

In practice we however do not know all the voltage and current density $(g,j)$ on the boundary $\partial\Omega$, we measure them at some discrete electrodes posed on a relatively open subset $\Gamma$ of the boundary only. An interpolation process is then required to derive the measured voltage and current density $(g_\delta,j_\delta)$ on $\Gamma$, where $\delta>0$ refers to the error level of the interpolation process and/or the measurements. The identification is now to reconstruct the electrical conductivity $q$ distributed inside the body $\Omega$ from boundary measurements of the voltage and current density, i.e. from the pair $(g_\delta,j_\delta)$. This problem is known to be
non-linear and severely ill-posed, due to the lack of data. To have an overview on numerical reconstruction for this identification problem one can examine, e.g., \cite{Knowles1998,KohnMcKenney90,LeRi08}.

In the present paper we investigate the problem of simultaneously identifying the conductivity (or diffusion) and the reaction (or absorption), subjecting several sets of measurement data on an accessible part of the boundary are available. The identification problem is a theme of diffuse optical tomography --- a category of applied sciences, e.g., neuroscience, medicine, wound monitoring and cancer detection. The interested reader may consult an incomplete list of references \cite{Ar99,ArSc09,Bo97,Ca07,Du10, GHA05,
He97,HeSo02,HoYa16,Ki96,LeYe13,Po95,RBH07} 
for detailed discussions on this issue.
In fact, 
assume $\Omega$ is an open, bounded and connected subset of $\mathbb{R}^d,~ d\ge 2$ with Lipschitz boundary $\partial\Omega$ and $\Gamma \subset \partial\Omega$ is a $(d-1)$-dimensional measurable surface. We consider the elliptic equation
\begin{align}
-\nabla \cdot \big(q \nabla u \big) + a u &= f \quad \mbox{~in} \quad \Omega,  \label{17-5-16ct1}\\
q \nabla u \cdot \vec{n} +\sigma u &= j^\dag \quad \mbox{on} \quad \Gamma,  \label{17-5-16ct1*}\\
q \nabla u \cdot \vec{n} +\sigma u &= j_0 \quad \mbox{on} \quad \partial\Omega\setminus\Gamma,  \label{17-5-16ct1**}\\
 u &= g^\dag \quad \mbox{on} \quad \Gamma,  \label{17-5-16ct1***}
\end{align}
where the source term $f\in H^{-1}(\Omega) := {H^1(\Omega)}^*$, the Neumann boundary condition $j_0\in H^{-1/2}(\partial\Omega\setminus\Gamma) := {\left(H^{1/2}(\partial\Omega\setminus\Gamma)\right)}^*$ (see \S \ref{25-5-20ct1} for the definition of Sobolev spaces on surfaces),  and the Robin coefficient $\sigma$ are assumed to be known with $\sigma \in L^\infty(\partial\Omega)$ and $\sigma(x)\ge 0$ a.e. on $\partial\Omega$.
The identification problem is to seek the pair $(q,a)$ in the aforementioned equation \eqref{17-5-16ct1} -- \eqref{17-5-16ct1***} assuming the full knowledge of 
all Cauchy data on $\Gamma$
\begin{align*}
\mathcal{C} := \left\{\big({q\nabla u\cdot \vec{n}}_{|\Gamma} + \sigma u_{|\Gamma}, u_{|\Gamma}\big) =: (j,g) \in H^{-1/2}(\Gamma) \times H^{1/2}(\Gamma)\right\}
\end{align*}
is given. Arridge and Lionheart  showed in \cite{ArLi05} the non-uniqueness of this identification problem for globally smooth identified coefficients. Nevertheless, Harrach in \cite{Ha09} (see also \cite{Ha12}) proved that the identification problem is uniquely solvable in the class of piecewise constant functions.

Our aim in this work is to reconstruct the pair $(q,a) \in Q\times A$ from several sets of measurement data $\big(j^i_\delta,g^i_\delta\big)_{i=1,\cdots,I} \subset H^{-1/2}(\Gamma) \times  H^{1/2}(\Gamma)$ of the exact data $\big(j^\dag,g^\dag\big)$ satisfying the noisy model
\begin{align}\label{31-5-17ct1}
\frac{1}{I} \left(\big\|j^i_\delta - j^\dag\big\|_{H^{-1/2}(\Gamma)} + \big\|g^i_\delta - g^\dag\big\|_{H^{1/2}(\Gamma)} \right)\le \delta
\end{align}
with $\delta >0$ standing for the error level of the observations. Here the admissible sets are assumed to be constrained of the general type
\begin{align}\label{14-11-18ct1}
Q := \{q \in L^\infty(\Omega) ~|~ \underline{q} \le q(x) \le \overline{q} \en \mbox{a.e. in} \en \Omega\}
\end{align}
and 
\begin{align}\label{14-11-18ct1*}
A := \{a\in L^\infty(\Omega) ~|~ \underline{a} \le a(x) \le \overline{a} \en \mbox{a.e. in} \en \Omega\},
\end{align}
where the constants $\underline{q},~ \overline{q},~ \underline{a},~ \overline{a}$ are given with $0< \underline{q} \le \overline{q}$ and $0< \underline{a} \le \overline{a}$. Furthermore, for simplicity of exposition, hereafter we assume that $I=1$, i.e. only one Neumann-Dirichlet pair $(j_\delta,g_\delta)$ available. We also discuss the multiple measurement in Section \ref{numerical examples}.

With the pair $(j_\delta,g_\delta)$ at hand we examine the Neumann boundary value problem
\begin{align}\label{17-5-16ct2}
-\nabla \cdot \big(q \nabla u \big) + a u = f \quad \mbox{~in} \quad \Omega,  \quad
q \nabla u \cdot \vec{n} +\sigma u = j_\delta \quad \mbox{on} \quad \Gamma,  \quad
q \nabla u \cdot \vec{n} +\sigma u = j_0 \quad \mbox{on} \quad \partial\Omega\setminus\Gamma
\end{align}
as well as the mixed boundary value problem
\begin{align}\label{17-5-16ct3}
-\nabla \cdot \big(q \nabla v \big) + a v = f \quad \mbox{~in} \quad \Omega,  \quad
v = g_\delta \quad \mbox{on} \quad \Gamma,  \quad
q \nabla v \cdot \vec{n} +\sigma v = j_0 \quad \mbox{on} \quad \partial\Omega\setminus\Gamma 
\end{align}
whose weak solutions are denoted by $N_{j_\delta}(q,a)$ and $M_{g_\delta}(q,a)$, respectively. Based on the variational approach of Kohn and Vogelius in \cite{Kohn_Vogelius1,Kohn_Vogelius11,Kohn_Vogelius2}, we propose the non-negative misfit energy functional
\begin{align*}
J_\delta(q,a) &:= \int_\Omega q\left|\nabla\big(N_{j_\delta}(q,a) - M_{g_\delta}(q,a)\big)\right|^2 dx + \int_\Omega a\big(N_{j_\delta}(q,a) - M_{g_\delta}(q,a)\big)^2 dx \\
&~\quad + \int_{\partial\Omega} \sigma\big(N_{j_\delta}(q,a) - M_{g_\delta}(q,a)\big)^2 ds
\end{align*}
for the identification problem and consider its minimizers over $Q\times A$ as reconstructions. However, since the identification problem is ill-posed, we make use a regularization method to seek stable solutions. Furthermore, for interests in estimating piecewise constant coefficients we therefore utilize total variation regularization combining with the quadratic stabilizing term, i.e. we consider the minimization problem
$$
\min_{(q,a)\in Q_{ad} \times A_{ad}} \Upsilon_{\delta,\rho} (q,a) := J_\delta(q,a) + \rho R(q,a), \eqno \left(\mathcal{P}_{\delta,\rho}\right)
$$
with
\begin{align*}
R(q,a) &:= \int_\Omega | \nabla q| + \int_\Omega | \nabla a| + \frac{1}{2}\|q\|^2_{L^2(\Omega)} + \frac{1}{2}\|a\|^2_{L^2(\Omega)}
\end{align*}
and 
$$Q_{ad} := Q \cap BV(\Omega) \en \mbox{and} \en A_{ad} := A \cap BV(\Omega),$$
where
$BV(\Omega)$ is the space of all functions of bounded total variation with the semi-norm $\int_\Omega |\nabla (\cdot)|$ and the norm $\int_\Omega |\nabla (\cdot)| + \| \cdot \|_{L^1(\Omega)}$ (cf.\ \S\ref{2-3-20ct1}), and $\rho >0$ is the regularization parameter. Total variation regularization originally introduced in image
denoising \cite{Rudin_Osher_Fatemi}. Somewhat later, it has been used to treat several
ill-posed and inverse problems over the
last decades, where the possibility of discontinuity in their solutions is interested in particular (cf.\ \cite{BurgerOsher13}).
We would like to mention that for coefficient identification problems in partial differential equations with smooth enough identified objects one may employ the regularization of Sobolev norms (see, e.g., \cite{aca93,cr14,JaGo12,KoLo88}). In the present paper we adopt the stabilized method of total variation combining with quadratic term first introduced in \cite{ChavenKu} for linear inverse problems 
to treat the non-linear identification problem with possibly discontinuous sought coefficients.

Let $V^h_1$ be the finite dimensional space of piecewise linear, continuous finite elements, and $N^h_{j_\delta} (q,a)$ and $M^h_{g_\delta}(q,a)$ be respectively the finite element approximations of $N_{j_\delta} (q,a)$ and $M_{g_\delta}(q,a)$ in $V^h_1$, where $h>0$ is the mesh size of the triangulation. We then approximate the problem $\left(\mathcal{P}_{\delta,\rho}\right)$ by the discrete one
$$\min_{(q,a)\in Q^h_{ad} \times A^h_{ad}} \Upsilon^h_{\delta,\rho} (q,a) :=  J^h_\delta(q,a) + \rho R(q,a), \eqno \big(\mathcal{P}^h_{\delta,\rho}\big)$$
where
\begin{align*}
J^h_\delta(q,a) &:= \int_\Omega q\left|\nabla\big(N^h_{j_\delta}(q,a) - M^h_{g_\delta}(q,a)\big)\right|^2 dx + \int_\Omega a\big(N^h_{j_\delta}(q,a) - M^h_{g_\delta}(q,a)\big)^2 dx\\
&~\quad + \int_{\partial\Omega} \sigma\big(N^h_{j_\delta}(q,a) - M^h_{g_\delta}(q,a)\big)^2 ds
\end{align*}
and 
$$Q^h_{ad} := Q \cap V^h_1 \subset Q \cap BV(\Omega) \en \mbox{and} \en A^h_{ad} := A \cap V^h_1 \subset A \cap BV(\Omega).$$
As the identification problem is non-linear and severely ill-posed, the stable analysis and convergence result of finite dimensional regularized solutions to the identification are crucial. 

Let the regularization parameter $\rho$ and the observation data $(j_\delta,g_\delta)$ be fixed and $\big(q^{h_n}_{\delta,\rho},a^{h_n}_{\delta,\rho}\big)$ denotes an arbitrary minimizer of $\big(\mathcal{P}^{h_n}_{\delta,\rho}\big)$ for each $n\in\mathbb{N}$, where $h_n \to 0$ as $n \to \infty$. We then show that the sequence $\big(q^{h_n}_{\delta,\rho},a^{h_n}_{\delta,\rho}\big)$ has a subsequence converging to an element $\big(q_{\delta,\rho},a_{\delta,\rho}\big) \in Q_{ad}\times A_{ad}$ in the $L^s(\Omega)^2$-norm for all $s \in [1,\infty)$ with $\big(q_{\delta,\rho},a_{\delta,\rho}\big)$ a solution of $\big(\mathcal{P}_{\delta,\rho}\big)$. 

Furthermore, let $\left(h_n\right)$ and  $\left(\delta_n\right)$ be any positive sequences converging to zero together with $\rho_n = \rho_n(h_n,\delta_n)$ being suitably chosen.
Assume that $\big(j_{\delta_n}, g_{\delta_n}\big) $ is a sequence satisfying
$$\big\|j_{\delta_n} - j^\dag\big\|_{H^{-1/2}(\Gamma)} + \big\|g_{\delta_n} - g^\dag\big\|_{H^{1/2}(\Gamma)} \le \delta_n$$
and that $\big(q^{h_n}_{\delta_n,\rho_n},a^{h_n}_{\delta_n,\rho_n}\big)$ is an arbitrary minimizer to $\big( \mathcal{P}_{\rho_n,\delta_n}^{h_n} \big)$ for each $n\in\N$. Then,

(i) There exist a subsequence of $\big(q^{h_n}_{\delta_n,\rho_n},a^{h_n}_{\delta_n,\rho_n}\big)$ denoted by the same symbol and a solution $(q^\dag,a^\dag)$ of the identification problem
$$
 \min_{\left\{(q,a)\in Q_{ad}\times A_{ad} ~\big|~ N_{j^\dag}(q,a) = M_{g^\dag}(q,a) \right\}} R(q,a) \eqno\left(\mathcal{IP}\right)
$$
such that $\big(q^{h_n}_{\delta_n,\rho_n},a^{h_n}_{\delta_n,\rho_n}\big)$ converges to $(q^\dag,a^\dag)$ in the $L^s(\Omega)^2$-norm for all $s \in [1,\infty)$ and
\begin{align*}
&\limn  \int_\Omega \big|\nabla q^{h_n}_{\delta_n,\rho_n}\big|  = \int_\Omega |\nabla q^\dag|  \quad\mbox{and}\quad \limn  \int_\Omega \big|\nabla a^{h_n}_{\delta_n,\rho_n}\big|  = \int_\Omega |\nabla a^\dag|,\\
&\limn D_{TV}^{\ell} \big(q^{h_n}_{\delta_n,\rho_n},q^\dag \big) = \limn D_{TV}^{\kappa} \big(a^{h_n}_{\delta_n,\rho_n},a^\dag \big) = 0
\end{align*}
for all
$(\ell,\kappa) \in \partial\left(\int_\Omega|\nabla(\cdot)|\right)(q^\dag) \times \partial\left(\int_\Omega|\nabla(\cdot)|\right)(a^\dag)$. Here $\partial\left(\int_\Omega|\nabla(\cdot)|\right)(\varphi)$ is the sub-differential of the semi-norm $\int_\Omega |\nabla(\cdot)|$ of the space $BV(\Omega)$ at $\varphi \in BV(\Omega)$ and $
D^{\ell}_{TV}(p,q)
$
is the Bregman distance with respect to $\int_\Omega |\nabla(\cdot)|$ and $\ell$ of
two elements $p, q$ (cf.\ \S\ref{2-3-20ct1}).

(ii) The sequences $\big( N^{h_n}_{j_{\delta_n}}\big(q^{h_n}_{\delta_n,\rho_n},a^{h_n}_{\delta_n,\rho_n}\big)\big) $ and $\big( M^{h_n}_{g_{\delta_n}}\big(q^{h_n}_{\delta_n,\rho_n},a^{h_n}_{\delta_n,\rho_n}\big)\big) $ converge in the $H^1(\Omega)$-norm to the solution $u(q^\dag,a^\dag)$. If the solution $(q^\dag,a^\dag)$ is uniquely defined, then the above convergences hold true for the whole sequence.

Furthermore, we show that the misfit term $J^h_\delta(q,a)$ is Fr\'echet differentiable and  for each $(q,a)\in Q^h_{ad}\times A^h_{ad}$, the Fr\'echet differential in the direction $(\eta_q,\eta_a) \in V^h_1 \times V^h_1$ given by
\begin{align*}
{J^h_\delta}'(q,a)(\eta_q,\eta_a) =  \int_\Omega \eta_q \left( \left| \nabla {M^h_{g_\delta}}(q,a) \right|^2 - \left| \nabla {N^h_{j_\delta}}(q,a) \right|^2 \right)dx +  \int_\Omega \eta_a \left( \left| {M^h_{g_\delta}}(q,a) \right|^2 - \left| {N^h_{j_\delta}}(q,a) \right|^2 \right)dx.
\end{align*}
Based on this fact, we perform some numerical results for the simultaneous coefficient identification problem, which illustrate the  efficiency of the proposed variational method.

To complete this introduction, we wish to mention that the problem of identifying the sole coefficient has been extensively investigated, see \cite{Aless,ChTa03,ChanTai,ChKu02,Chenzou,
fal83,Haoq,hao_quyen1,hao_quyen2,HKQ18,ito-kunisch2008,KeungZou,Qu19,Ric,
Vainikko_Kunisch1993,zou1998} and many others in the literature. We have not yet found investigations for the multiple coefficient identification problem with {\it boundary observations}, however with {\it distributed observations} in \cite{Baku,hao_quyen3,hein}.  By using a non-standard version of the misfit term combining with an appropriate regularized technique we could in the present paper outline that two coefficients distributed inside the physical domain can be simultaneously reconstructed from a finite number of observations on a part of the boundary.

The paper is organized as follows. In Section \ref{12-9-19ct1} we introduce some useful notations and show the existence of a minimizer of the regularized minimization problem. Finite element method for the identification problem is presented in Section \ref{discrete}. Stability analysis of the proposed regularization approach and convergence of the finite dimensional approximations to the identification are enclosed in Section \ref{stability}. We in Section \ref{29-8-19ct3} perform the differentials of the discrete coefficient-to-solution operators and of the associated cost functional together with a projected gradient method to reach minimizers of the formulated identification problems. Finally, some numerical examples supporting our analytical findings are presented in Section \ref{numerical examples}. 

\section{Preliminaries}\label{12-9-19ct1}

\subsection{Sobolev spaces on surfaces} \label{25-5-20ct1}

Assume that $\Gamma \subset \partial\Omega$ is a $(d-1)$-dimensional measurable surface and $s \in (0,1)$, we denote by (see, \cite[\S 1.2.1.4]{Pechstein})
\begin{align*}
H^s(\Gamma) := \left\{ g\in L^2(\Gamma) ~|~ \|g\|_{H^s(\Gamma)} < \infty \right\}
\end{align*}
with the norm
\begin{align*}
\|g\|_{H^s(\Gamma)} := \left(\|g\|^2_{L^2(\Gamma)} + |g|^2_{H^s(\Gamma)} \right)^{1/2}
\end{align*}
and the semi-norm
\begin{align*}
|g|_{H^s(\Gamma)} := \left(\int_\Gamma \int_\Gamma \frac{|g(x) - g(y)|^2}{|x-y|^{d-1+2s}} dx dy\right)^{1/2}.
\end{align*}
We mention that, equipped with the norm $\|\cdot\|_{H^s(\Gamma)}$, $H^s(\Gamma)$ is a Hilbert space for all $s \in (0,1)$. Further, one may consult \cite[\S 2.4]{SaSc11} for the definition of the space $H^s(\Gamma)$ with $s \ge 1$. 

We also note that (see, \cite[Lemma 1.21]{Pechstein}) if $\Gamma$ is open, then
\begin{align*}
H^s(\Gamma) = \left\{ g\in L^2(\Gamma) ~|~ \exists G \in H^s(\partial\Omega),~ G_{|\Gamma} = g \right\}
\end{align*}
and the norm
\begin{align*}
\inf_{G \in H^s(\partial\Omega),~ G_{|\Gamma} = g} \|G\|_{H^s(\partial\Omega)}
\end{align*}
is equivalent to the norm $\|\cdot\|_{H^s(\Gamma)}$.

Finally, let us denote by
$$\gamma : H^1(\Omega) \to H^{1/2}(\partial\Omega)$$
the continuous Dirichlet trace operator with
$$\gamma^{-1} : H^{1/2}(\partial\Omega) \to H^1(\Omega)$$
its continuous right inverse, i.e. $ (\gamma\circ \gamma^{-1}) g =g$ for all $g\in H^{1/2}(\partial\Omega)$. Moreover, if $\Omega$ is of the $C^{k-1,1}$ class for some $k\ge 1$, then $\gamma : H^s(\Omega) \to H^{s-1/2}(\partial\Omega)$ is continuous for all $s\in (0,k]$ (see, e.g., \cite[Theorem 1.23]{Pechstein}).

\subsection{Bregman distance with respect to total variation semi-norm}\label{2-3-20ct1}

We start with briefly summarizing the work
\cite{Resmer} about the Bregman
distance related to a proper convex function. Let $\mathcal{B}$ be
a Banach space with $\mathcal{B}^*$ the dual space and $\Phi:
\mathcal{B}\rightarrow (-\infty, +\infty]$ is a proper convex
function, i.e. $\text{Dom}{\Phi} :=\{ q\in\mathcal{B}~|~\Phi(q)<+\infty\}
\neq\emptyset$.

Let $\partial \Phi(q)$ stand for the sub-differential of $\Phi$
at $q\in\text{Dom}{\Phi}$ defined by
$$
\partial \Phi(q) := \{q^*\in \mathcal{B}^*~|~
\Phi(p)\ge \Phi(q)+\<q^*, p-q\>_{(\mathcal{B}^*, \mathcal{B})}  \text{~for
all~} p\in \mathcal{B}\}.
$$
The set $\partial \Phi(q)$ may be empty; however, if $\Phi$ is  continuous
at $q$, then it is non-empty. Further, $\partial \Phi(q)$ is convex and
weak* compact (see, \cite{eke}).
In case $\partial \Phi(q) \neq\emptyset$, we for any fixed $p\in
\mathcal{B}$ denote by
$$
D_\Phi (p,q) :=\left\{\Phi(p) - \Phi(q)  - \<q^*, p-q\>_{(\mathcal{B}^*,
\mathcal{B})} ~\big|~ q^*\in \partial \Phi(q)\right\}
$$
and with a fixed element $q^*\in \partial \Phi(q)$ the non-negative quantity
$$
D^{q^*}_\Phi (p,q) := \Phi(p) - \Phi(q) - \<q^*, p-q\>_{(\mathcal{B}^*,
\mathcal{B})}
$$
is called the Bregman distance with respect to $\Phi$ and $q^*$ of
two elements $p, q\in \mathcal{B}$.

The Bregman distance is not a metric on $\mathcal{B}$ in general.
However, $D^{q^*}_\Phi (q,q)=0$, and in case
$\Phi$ is the strictly convex function, the identity $D^{q^*}_\Phi (p,q)= 0$ implies $p=q$. The notion of Bregman distance was first given by
Bregman \cite{Bregman} for Fr\'echet differentiable $\Phi$ and it was
generalized by Kiwiel \cite{Kiw} to nonsmooth but strictly convex
$\Phi$. Burger and Osher \cite{BurOs}  further generalized this notion
for $\Phi$ being neither smooth, nor strictly convex.

Next, we present the definition of functions with bounded total variation;  for more details, we may consult \cite{Evan,Giusti}. A function $q\in
L^{1}(\Omega)$ is said to be of bounded total variation if
\begin{align*}
TV(q) := \int_{\Omega}\left|\nabla q\right|
:=\sup\left\{\int_{\Omega} q~\divv~gdx~|~g\in C^1_c(\Omega)^d,~
\left|g(x)\right|_{\infty} \le
1,~x\in\Omega\right\}<\infty,
\end{align*}
where $\left|\cdot\right|_{\infty}$ denotes the $\ell_{\infty}$-norm
on $\mathbb{R}^d$, i.e. 
$\left|x\right|_{\infty}=\max\limits_{1\le i\le d}\left|x_i\right|$ for all $x= (x_1,\ldots,x_d) \in \mathbb{R}^d$.
The space of all functions in $L^{1}(\Omega)$ with bounded total
variation  is denoted by
$$BV(\Omega)=\left\{q\in L^{1}(\Omega) ~\big|~
\int_{\Omega}\left|\nabla q\right|<\infty\right\}.$$
It is a Banach space endowed with the norm
$$
\left\|q\right\|_{BV(\Omega)}
:=\left\|q\right\|_{L^{1}(\Omega)}+\int_{\Omega}\left|\nabla
q\right|,
$$
while $\int_{\Omega}\left|\nabla
(\cdot)\right|$ is a semi-norm of $BV(\Omega)$.
Furthermore, if
$\Omega$ is an open bounded set with
Lipschitz boundary, then $W^{1,1}(\Omega)\varsubsetneq BV(\Omega)$.

Let $\partial\big(\int_\Omega |\nabla(\cdot)|\big)(q)$ be the sub-differential of the semi-norm $\int_\Omega |\nabla(\cdot)|$ at $q \in BV(\Omega)$. As $q \mapsto \int_\Omega |\nabla q|$ is a continuous functional on the space $BV(\Omega)$, the set $\partial\big(\int_\Omega |\nabla(\cdot)|\big)(q) \neq \emptyset$.
Then for a fixed element $q^*\in \partial\big(\int_\Omega |\nabla(\cdot)|\big)(q)$ the Bregman distance with respect to $\int_\Omega |\nabla(\cdot)|$ and $q^*$ of
two elements $p, q$ reads as
$$
D^{q^*}_{TV}(p,q) := \int_\Omega |\nabla p| - \int_\Omega |\nabla q| - \<q^*, p -q\>_{\big({BV(\Omega)}^*,BV(\Omega)\big)}.
$$

\subsection{Auxiliary results}

The expression 
\begin{align*}
\int_\Omega q \nabla u \cdot \nabla vdx + \int_\Omega auvdx + \int_{\partial\Omega} \sigma uvds
\end{align*}
generates an inner product on the space
$ H^1(\Omega)$ which is equivalent to the usual one, i.e. there exist positive constants $c_1, c_2$ such that
\begin{align}\label{18-10-16ct1}
c_{1}\|u\|^2_{H^1(\Omega)} \le \int_\Omega q \nabla u \cdot \nabla udx + \int_\Omega au^2dx + \int_{\partial\Omega} \sigma u^2ds \le c_2\|u\|^2_{H^1(\Omega)}
\end{align}
for all $u\in H^1(\Omega)$. 
Therefore, for each $(q,a) \in Q\times A$ the Neumann boundary value problem \eqref{17-5-16ct2}
defines a unique weak solution denoted by $N_{j_\delta}(q,a)$ in the sense that $N_{j_\delta}(q,a)\in H^1(\Omega)$ and the equation
\begin{align}\label{17-10-16ct2}
\int_\Omega q\nabla N_{j_\delta}(q,a) \cdot \nabla \phi dx  
&+ \int_\Omega a N_{j_\delta}(q,a)  \phi dx + \int_{\partial\Omega} \sigma N_{j_\delta}(q,a) \phi ds \\
&= \<f,\phi\>_{\big(H^{-1}(\Omega), H^1(\Omega) \big)} + \<j_\delta, \phi\>_{\big( H^{-1/2}(\Gamma), H^{1/2}(\Gamma) \big)} + \<j_0, \phi\>_{\big( H^{-1/2}(\partial\Omega\setminus\Gamma), H^{1/2}(\partial\Omega\setminus\Gamma) \big)} \notag
\end{align}
is satisfied for all $\phi\in H^1(\Omega)$. Furthermore, there holds the estimate
\begin{align}\label{17-10-16ct4}
\|N_{j_\delta}(q,a)\|_{H^1(\Omega)} \le C \left( \|f\|_{H^{-1}(\Omega)}  + \|j_\delta\|_{H^{-1/2}(\Gamma)} + \|j_0\|_{H^{-1/2}(\partial\Omega\setminus\Gamma)}\right)
\end{align}
for some positive constant $C$. A function $M_{g_\delta}(q,a)$ is said to be a unique weak solution of the mixed boundary value problem \eqref{17-5-16ct3}
if $M_{g_\delta}(q,a)\in H^1(\Omega)$ with $M_{g_\delta}(q,a)_{|\Gamma} = g_\delta$ and the equation
\begin{align}\label{17-10-16ct2*}
\int_\Omega q\nabla M_{g_\delta}(q,a) \cdot \nabla \phi dx  
&+ \int_\Omega a M_{g_\delta}(q,a)  \phi dx + \int_{\partial\Omega} \sigma M_{g_\delta}(q,a) \phi ds \\
&= \<f,\phi\>_{\big(H^{-1}(\Omega), H^1(\Omega) \big)} + \<j_0, \phi\>_{\big( H^{-1/2}(\partial\Omega\setminus\Gamma), H^{1/2}(\partial\Omega\setminus\Gamma) \big)} \notag
\end{align}
is fulfilled for all $\phi\in H^1_0(\Omega\cup\Gamma)$, where 
$$
H^1_0(\Omega\cup\Gamma) := \overline{C^\infty_c(\Omega\cup\Gamma)}^{H^1(\Omega)} = \{\phi \in H^1(\Omega) ~|~ \phi_{|\Gamma} = 0\}$$
the bar denoting the closure in $H^1(\Omega)$ and $C^\infty_c(\Omega\cup\Gamma)$ consisting all functions $\phi\in C^\infty(\Omega)$ with $\mbox{supp}\phi$ being a compact subset of $\Omega\cup\Gamma$ (cf.\ \cite{Troianiello}). The above weak solution satisfies  the estimate
\begin{align}\label{14-11-18ct3}
\|M_{g_\delta}(q,a)\|_{H^1(\Omega)} \le  C\left(\|f\|_{H^{-1}(\Omega)} + \|g_\delta\|_{H^{1/2}(\Gamma)}  + \|j_0\|_{H^{-1/2}(\partial\Omega\setminus\Gamma)}\right).
\end{align}

Thus we define the {\it non-linear} coefficient-to-solution operators
$$N_{j_\delta},~ M_{g_\delta}:~ Q\times A \subset L^\infty(\Omega) \times L^\infty(\Omega) \to H^1(\Omega)$$
which are uniformly bounded, due to \eqref{17-10-16ct4} and \eqref{14-11-18ct3}.

We here present some properties of the coefficient-to-solution operators.

\begin{lemma}\label{weakly conv.}
Assume that the sequence $\left( q_n,a_n\right)\subset Q\times A$
converges to $(q,a)$ almost everywhere in $\Omega$. Then $(q,a)\in Q\times A$ and the sequence $\big(N_{j_\delta}(q_n,a_n),M_{g_\delta}(q_n,a_n)\big)$ converges to $\big(N_{j_\delta}(q,a),M_{g_\delta}(q,a)\big)$ in the $H^1(\Omega)\times H^1(\Omega)$-norm.
\end{lemma}

\begin{proof}
By the equation \eqref{17-10-16ct2},  we for each $n\in \mathbb{N}$ get that
\begin{align*}
&\int_\Omega q_n\nabla \big(N_{j_\delta}(q_n,a_n) - N_{j_\delta}(q,a)\big) \cdot \nabla \phi dx  
+ \int_\Omega a_n \big(N_{j_\delta}(q_n,a_n) - N_{j_\delta}(q,a)\big)  \phi dx \\
&~\quad + \int_{\partial\Omega} \sigma \big(N_{j_\delta}(q_n,a_n) - N_{j_\delta}(q,a)\big) \phi ds
= \int_\Omega (q-q_n)\nabla N_{j_\delta}(q,a) \cdot \nabla \phi dx + \int_\Omega (a - a_n) N_{j_\delta}(q,a) \phi dx.
\end{align*}
Taking $\phi = N_{j_\delta}(q_n,a_n) - N_{j_\delta}(q,a)$ and using the inequality \eqref{18-10-16ct1}, we arrive at
\begin{align*}
\|N_{j_\delta}(q_n,a_n) - N_{j_\delta}(q,a)\|_{H^1(\Omega)} 
\le C \left( \sqrt{\int_\Omega |q - q_n|^2 |\nabla N_{j_\delta}(q,a)|^2 dx} + \sqrt{\int_\Omega |a - a_n|^2 |N_{j_\delta}(q,a)|^2 dx} \right)
\end{align*}
By the Lebesgue dominated convergence theorem, we conclude that $\limn \|N_{j_\delta}(q_n,a_n) - N_{j_\delta}(q,a)\|_{H^1(\Omega)} =0$. Similarly, we also obtain $\limn \|M_{g_\delta}(q_n,a_n) - M_{g_\delta}(q,a)\|_{H^1(\Omega)} =0$, which finishes the proof.
\end{proof}

Next, let us quote the following useful results.

\begin{lemma}[\cite{attouch}]\label{bv1}

(i) Let $\left(w_n\right)$ be a bounded sequence in the $BV(\Omega)$-norm.
Then a subsequence not relabeled and an element
$w\in BV(\Omega)$ exist such that $\left(w_n\right)$ converges to $w$ in the
$L^1(\Omega)$-norm.

(ii) Let $\left(w_n\right)$ be a sequence in $BV(\Omega)$ converging to
$w$ in  the $L^1(\Omega)$-norm. Then $w \in BV(\Omega)$ and
\begin{align*}
\int_\Omega \left|\nabla w\right| \le \liminfn \int_\Omega|\nabla w_n|.
\end{align*}
\end{lemma}

\begin{lemma}[\cite{bartels}]\label{appr.BV}
Assume that $w\in BV(\Omega)$. Then for all $\epsilon >0$ an element $w^\epsilon \in C^\infty(\Omega)$ exists such that
$$\int_\Omega|w-w^\epsilon| \le \epsilon\int_\Omega|\nabla w|,~ \int_\Omega|\nabla w^\epsilon| \le (1+C\epsilon)\int_\Omega|\nabla w| \mbox{~and~} \int_\Omega|D^2 w^\epsilon| \le C\epsilon^{-1}\int_\Omega|\nabla w|,$$
where the positive constant $C$ is independent of $\epsilon$.
\end{lemma}

We are now in the position to prove the main result of the section.

\begin{theorem}\label{existance}
The problem $\left(\mathcal{P}_{\delta,\rho}\right)$ attains a solution $\big(q_{\delta,\rho},a_{\delta,\rho}\big)$, which is called the regularized solution of the identification problem.
\end{theorem}

\begin{proof}
Let $\left(q_n,a_n\right) \subset Q_{ad} \times A_{ad}$ be a minimizing sequence of the problem $\left(\mathcal{P}_{\delta,\rho}\right)$, i.e.
\begin{align*}
\limn \Upsilon_{\delta,\rho} (q_n,a_n) 
= \inf_{(q,a)\in Q_{ad} \times A_{ad}} \Upsilon_{\delta,\rho} (q,a).
\end{align*}
Therefore, the sequence $\left(q_n,a_n\right)$ is bounded in the $BV(\Omega)$-norm.
By Lemma \ref{bv1}, a subsequence which is not relabeled and an element $(q,a) \in Q_{ad} \times A_{ad}$ exist such that
\begin{quote}
$\left(q_n,a_n\right)$ converges to $(q,a)$ in the
$L^1(\Omega)\times L^1(\Omega)$-norm, \\
$\left(q_n,a_n\right)$ converges to $(q,a)$ almost everywhere in $\Omega$,\\
$\int_\Omega |\nabla q| \le \liminfn \int_\Omega |\nabla q_n|$ and $\int_\Omega |\nabla a| \le \liminfn \int_\Omega |\nabla a_n|$.
\end{quote}
By the inequality
\begin{align*}
\|q_n -q \|^2_{L^2(\Omega)} + \|a_n -a \|^2_{L^2(\Omega)}
\le 2\max(\overline{q} , \overline{a}) (\|q_n -q \|_{L^1(\Omega)} + \|a_n -a \|_{L^1(\Omega)}),
\end{align*}
the sequence $\left(q_n,a_n\right)$ also converges to $(q,a)$ in the
$L^2(\Omega)\times L^2(\Omega)$-norm. We thus have
\begin{align}\label{26-3-16ct5}
R(q,a) \le \liminfn R(q_n,a_n).
\end{align}
Furthermore, an application of Lemma \ref{weakly conv.} deduces that the sequence $\big(N_{j_\delta}(q_n,a_n),M_{g_\delta}(q_n,a_n)\big)$ converges to $\big(N_{j_\delta}(q,a),M_{g_\delta}(q,a)\big)$ in the $H^1(\Omega)\times H^1(\Omega)$-norm and then
\begin{align}\label{15-11-18ct1}
J_\delta(q,a) = \limn J_\delta(q_n,a_n).
\end{align}
Therefore, we obtain from \eqref{26-3-16ct5} -- \eqref{15-11-18ct1} that
\begin{align*}
\Upsilon_{\delta,\rho} (q,a)
&\le \limn J_\delta\left(q_n,a_n\right) + \liminfn \rho R(q_n,a_n) \\
&= \liminfn \left(J_\delta\left(q_n,a_n\right) + \rho R(q_n,a_n)\right) \\
&= \inf_{(q,a)\in Q_{ad} \times A_{ad}} \Upsilon_{\delta,\rho} (q,a)
\end{align*}
and $(q,a)$ is hence a solution of the problem $\left(\mathcal{P}_{\delta,\rho}\right)$, which finishes the proof.
\end{proof}

\section{Finite element discretization}\label{discrete}

Hereafter we assume that $\Omega$ is a Lipschitz polygonal domain and $\left(\mathcal{T}^h\right)_{0<h<1}$ is a quasi-uniform family of regular triangulations of $\overline{\Omega}$ with the mesh size $h$ such that each vertex of the polygonal boundary $\partial\Omega$ is a node of $\mathcal{T}^h$. Let us denote by
\begin{align*}
&V_1^h := \left\{v^h\in C(\overline\Omega)
~|~{v^h}_{|T} \in \mathcal{P}_1, ~~\forall
T\in \mathcal{T}^h\right\}\\
&V_{1,0}^h := V_1^h \cap H^1_0(\Omega\cup\Gamma),
\end{align*}
where $\mathcal{P}_{1}$ consists of all polynomial functions of degree
less than or equal to $1$. For each $(q,a)\in Q\times A$ the variational equations
\begin{align}\label{10/4:ct1}
\int_\Omega q\nabla u^h \cdot \nabla \phi^h dx  
&+ \int_\Omega a u^h  \phi^h dx + \int_{\partial\Omega} \sigma u^h \phi^h ds \\
&= \<f,\phi^h\>_{\big(H^{-1}(\Omega), H^1(\Omega) \big)} + \<j_\delta, \phi^h\>_{\big( H^{-1/2}(\Gamma), H^{1/2}(\Gamma) \big)} + \<j_0, \phi^h\>_{\big( H^{-1/2}(\partial\Omega\setminus\Gamma), H^{1/2}(\partial\Omega\setminus\Gamma) \big)} \notag
\end{align}
for all $\phi^h\in V^h_1$ and
\begin{align}\label{10/4:ct1*}
\int_\Omega q\nabla v^h \cdot \nabla \phi^h dx  
&+ \int_\Omega a v^h  \phi^h dx + \int_{\partial\Omega} \sigma v^h \phi^h ds \\
&= \<f,\phi^h\>_{\big(H^{-1}(\Omega), H^1(\Omega) \big)}  + \<j_0, \phi^h\>_{\big( H^{-1/2}(\partial\Omega\setminus\Gamma), H^{1/2}(\partial\Omega\setminus\Gamma) \big)} \notag
\end{align}
for all $\phi^h\in V_{1,0}^h$ and ${v^h}_{|\Gamma} = g_\delta$
admit unique solutions $u^h := N^h_{j_\delta}(q,a) \in V_1^h$ and $v^h := M^h_{g_\delta}(q,a) \in V_1^h$, respectively. Furthermore, the estimates
\begin{align}
\|N^h_{j_\delta}(q,a)\|_{H^1(\Omega)} 
&\le C \left( \|f\|_{H^{-1}(\Omega)}  + \|j_\delta\|_{H^{-1/2}(\Gamma)} + \|j_0\|_{H^{-1/2}(\partial\Omega\setminus\Gamma)}\right)  \label{18/5:ct1}\\
\|M^h_{g_\delta}(q,a)\|_{H^1(\Omega)} 
&\le C \left( \|f\|_{H^{-1}(\Omega)}  + \|g_\delta\|_{H^{1/2}(\Gamma)} + \|j_0\|_{H^{-1/2}(\partial\Omega\setminus\Gamma)}\right)  \label{18/5:ct1*}
\end{align}
hold true, where the positive constant $C$ is independent of $h$.

\begin{remark}\label{27-12-18ct2}
Due to the standard theory of the finite element method for elliptic problems (cf.\ \cite{Brenner_Scott}), we for any fixed $(q,a)\in Q\times A$ get the limits 
\begin{align}\label{27-9-19ct1}
\lim_{h\to 0} \big\|N^h_{j_\delta}(q,a) - N_{j_\delta}(q,a)\big\|_{H^1(\Omega)} = \lim_{h\to 0} \big\|M^h_{g_\delta}(q,a) - M_{g_\delta}(q,a)\big\|_{H^1(\Omega)} =0.
\end{align}
Furthermore, under additional assumptions $q \in {C^{0,1}(\overline{\Omega})}$, $a \in L^\infty(\Omega)$, $f \in L^2(\Omega)$, $\sigma\in L^\infty(\partial\Omega)$, $j_0 \in H^{1/2}(\partial\Omega\setminus\Gamma)$, $j_\delta \in H^{1/2}(\Gamma)$,  $g_\delta \in H^{3/2}(\Gamma)$ and $\Omega$ is either of the class $C^{1,1}$ with the open portion $\Gamma$ being also closed (\cite[Theorem 2.24]{Troianiello}) or a Lipschitz polygonal domain (\cite[Theorem 4.3.1.4]{Grisvad}, see also \cite[Theorem 3.2.5]{BaRo89}), the weak solutions $N_{j_\delta}(q,a),~ M_{g_\delta}(q,a)\in H^2(\Omega)$ satisfying
\begin{align*}
\|N_{j_\delta}(q,a)\|_{H^2(\Omega)} 
&\le C \left( \|f\|_{L^2(\Omega)}  + \|j_\delta\|_{H^{1/2}(\Gamma)} + \|j_0\|_{H^{1/2}(\partial\Omega\setminus\Gamma)}\right)  \\
\|M_{g_\delta}(q,a)\|_{H^2(\Omega)} 
&\le C \left( \|f\|_{L^2(\Omega)}  + \|g_\delta\|_{H^{3/2}(\Gamma)} + \|j_0\|_{H^{1/2}(\partial\Omega\setminus\Gamma)}\right) 
\end{align*}
which yield the error bounds
\begin{align}
\big\|N^h_{j_\delta}(q,a) - N_{j_\delta}(q,a)\big\|_{L^2(\Omega)} + h\big\|N^h_{j_\delta}(q,a) - N_{j_\delta}(q,a)\big\|_{H^1(\Omega)} 
& \le Ch^2\|N_{j_\delta}(q,a)\|_{H^2(\Omega)} \label{8-11-18ct1}\\
\big\|M^h_{g_\delta}(q,a) - M_{g_\delta}(q,a)\big\|_{L^2(\Omega)} + h\big\|M^h_{g_\delta}(q,a) - M_{g_\delta}(q,a)\big\|_{H^1(\Omega)} 
&\le Ch^2\|M_{g_\delta}(q,a)\|_{H^2(\Omega)}. \label{8-11-18ct2}
\end{align}
\end{remark}

We introduce the Lagrange nodal value interpolation operator 
$$I^h_1: C(\overline{\Omega}) \to V^h_1.$$
By the continuous embedding $W^{1,p}(\Omega)\hookrightarrow C(\overline{\Omega})$ with $p>d$, the operator $I^h_1 : W^{1,p}(\Omega)\to V^h_1$ is well defined. Furthermore, see, e.g., \cite{Brenner_Scott}, it holds the limit 
\begin{align}\label{16-11-18ct1}
\lim_{h\to 0}\|I^h_1\phi - \phi\|_{W^{1,p}(\Omega)} =0
\end{align}
and the estimate
\begin{align}\label{16-11-18ct2}
\|I^h_1\phi - \phi\|_{L^p(\Omega)} \le  Ch\|\phi\|_{W^{1,p}(\Omega)}.
\end{align}

We have the following existence result. Its proof exactly follows  as in the continuous case, is therefore omitted here.
\begin{theorem}\label{dis-existance}
The discrete regularized problem 
$\big(\mathcal{P}^h_{\delta,\rho}\big)$
attains a minimizer $\big(q^h_{\delta,\rho}, a^h_{\delta,\rho}\big)$, which is called the discrete regularized solution of the identification problem.
\end{theorem}

\section{Convergence analysis}\label{stability}

The aim of this section is to prove the stability of the proposed regularization approach and the convergence of finite element approximations to the identification.

\begin{theorem}\label{29-8-19ct5}
Assume that the regularization parameter $\rho$ and the observation data $(j_\delta,g_\delta)$ are fixed. For each $n\in\mathbb{N}$ let $(q_n,a_n) := \big(q^{h_n}_{\delta,\rho},a^{h_n}_{\delta,\rho}\big)$ denote an arbitrary minimizer of $\big(\mathcal{P}^{h_n}_{\delta,\rho}\big)$, where $h_n \to 0$ as $n\to \infty$. Then the sequence $(q_n,a_n)$ has a subsequence not relabeled converging to an element $\big(q_{\delta,\rho},a_{\delta,\rho}\big) \in Q_{ad}\times A_{ad}$ in the $L^s(\Omega)$-norm for all $s \in [1,\infty)$. Furthermore, 
\begin{align}
&\limn  \int_\Omega |\nabla q_n|  = \int_\Omega |\nabla q_{\delta,\rho}|  \quad\mbox{and}\quad \limn  \int_\Omega |\nabla a_n|  = \int_\Omega |\nabla a_{\delta,\rho}|,\label{28-8-19ct2}\\
&\limn D_{TV}^{\ell}(q_n,q_{\delta,\rho}) = \limn D_{TV}^{\kappa}(a_n,a_{\delta,\rho}) = 0 \label{29-8-19ct4}
\end{align}
for all
$(\ell,\kappa) \in \partial\left(\int_\Omega|\nabla(\cdot)|\right)(q_{\delta,\rho}) \times \partial\left(\int_\Omega|\nabla(\cdot)|\right)(a_{\delta,\rho})$, where
$\big(q_{\delta,\rho},a_{\delta,\rho}\big)$ is a minimizer of $\big(\mathcal{P}_{\delta,\rho}\big)$.
\end{theorem}

\begin{proof}
Let $(q,a) \in Q_{ad}\times A_{ad}$ be arbitrary but fixed. Due to Lemma \ref{appr.BV}, for any fixed $\epsilon\in(0,1)$ an element  $\big(q^\epsilon, a^\epsilon\big) \in C^\infty(\Omega) \times C^\infty(\Omega)$ exists such that
\begin{equation}\label{28-3-18ct1}
\begin{aligned}
\frac{1}{2}\left\|q-q^\epsilon\right\|^2_{L^2(\Omega)} + \frac{1}{2}\left\|a -a^\epsilon\right\|^2_{L^2(\Omega)} 
&\le \epsilon\max(\overline{q},\overline{a}) \left( \int_\Omega |\nabla q| + \int_\Omega |\nabla a| \right) \\
\int_\Omega |\nabla q^\epsilon| + \int_\Omega |\nabla a^\epsilon|
&\le (1+C\epsilon)\left( \int_\Omega |\nabla q| + \int_\Omega |\nabla a| \right)
\end{aligned}
\end{equation}
for some positive constant $C$ independent of $\epsilon$. We denote by
$$r^\epsilon (q) := r^\epsilon := \max\left( \underline{q}, \min\left(q^\epsilon, \overline{q} \right) \right) \en \mbox{and} \en b^\epsilon(a) := b^\epsilon := \max\left( \underline{a}, \min\left(a^\epsilon, \overline{a} \right) \right)$$
that satisfy that 
$$\big(r^\epsilon, b^\epsilon\big) \in \big(Q \cap W^{1,\infty}(\Omega)\big) \times \big(A \cap W^{1,\infty}(\Omega)\big)$$
and
$$\big(r^\epsilon_n, b^\epsilon_n\big) := \big(I^{h_n}_1 r^\epsilon, I^{h_n}_1 b^\epsilon\big)  \in Q^h_{ad} \times A^h_{ad}.$$
Let $p>d$ and $p^*$ be the adjoint number of $p$, i.e. $\frac{1}{p} + \frac{1}{p^*} = 1$. We get
\begin{align}\label{6-4-16ct5}
|\Omega|^{-1/p^*} \int_\Omega |\nabla r^\epsilon_n|
&= |\Omega|^{-1/p^*} \int_\Omega |\nabla r^\epsilon_n(x)|dx \le \left( \int_\Omega |\nabla r^\epsilon_n(x)|^pdx \right)^{1/p}=  |r^\epsilon_n|_{W^{1,p}(\Omega)} \notag\\
&\le C + |r^\epsilon|_{W^{1,p}(\Omega)}
\end{align}
for $n$ large enough, by the limit \eqref{16-11-18ct1}. Furthermore, using \eqref{28-3-18ct1}, we have the estimate
\begin{align}\label{6-4-16ct5*}
|r^\epsilon|_{W^{1,p}(\Omega)} 
&= \int_\Omega |\nabla r^\epsilon(x)|dx
= \int_{ \{x\in\Omega ~|~ r^\epsilon(x) = q^\epsilon(x)\}} |\nabla r^\epsilon(x)|dx \le \int_\Omega |\nabla q^\epsilon(x)|dx \notag\\
&\le (1+C\epsilon)\left( \int_\Omega |\nabla q| + \int_\Omega |\nabla a| \right),
\end{align}
by the fact that $r^\epsilon$ is constant on $\{x\in\Omega ~|~ r^\epsilon(x) \not= q^\epsilon(x)\}$. Combining \eqref{6-4-16ct5} and \eqref{6-4-16ct5*}, we have the boundedness 
\begin{align}\label{16-11-18ct3}
\int_\Omega |\nabla r^\epsilon_n| + \int_\Omega |\nabla b^\epsilon_n| \le C
\end{align}
for all $n\in\mathbb{N}$ and $\epsilon\in (0,1)$. 

Now, by the definition of $(q_n,a_n)$, we for all $n\in \mathbb{N}$ get  that
\begin{align}\label{29-8-19ct1}
J^{h_n}_{\delta} \left(q_n,a_n\right) + \rho R(q_n,a_n)
\le J^{h_n}_{\delta} \left(r^\epsilon_n,b^\epsilon_n\right) + \rho R(r^\epsilon_n,b^\epsilon_n).
\end{align}
By \eqref{18/5:ct1} and \eqref{18/5:ct1*}, it holds 
$
J^{h_n}_{\delta} \left(r^\epsilon_n,b^\epsilon_n\right) \le C.
$
We thus deduce from \eqref{16-11-18ct3} -- \eqref{29-8-19ct1} that
$$R(q_n,a_n) \le C$$
for all $n\in \mathbb{N}$. 
An application of Lemma \ref{bv1} then follows that a subsequence of $(q_n,a_n)$ not relabeled and an element $(\widehat{q},\widehat{a})\in Q_{ad}\times A_{ad}$ exist such that $\left(q_n,a_n\right)$ converges to $(\widehat{q},\widehat{a})$ in the $L^1(\Omega)$-norm and
\begin{align*}
\int_\Omega \left|\nabla \widehat{q}\right| \le \liminfn \int_\Omega|\nabla q_n| \quad\mbox{and}\quad 
\int_\Omega \left|\nabla \widehat{a}\right| \le \liminfn \int_\Omega|\nabla a_n|.
\end{align*}
By the inequalities
\begin{align*}
\|q_n - \widehat{q}\|^s_{L^s(\Omega)} = \int_\Omega |q_n - \widehat{q}|\cdot|q_n - \widehat{q}|^{s-1} \le \int_\Omega |q_n - \widehat{q}|\cdot\left(|q_n|+|\widehat{q}|\right)^{s-1} \le \left(2\overline{q}\right)^{s-1}\|q_n - \widehat{q}\|_{L^1(\Omega)}
\end{align*}
and
\begin{align*}
\|a_n - \widehat{a}\|^s_{L^s(\Omega)} \le \left(2\overline{a}\right)^{s-1}\|a_n - \widehat{a}\|_{L^1(\Omega)},
\end{align*}
we deduce that $\left(q_n,a_n\right)$ in fact converges to $(\widehat{q},\widehat{a})$ in the $L^s(\Omega)$-norm for all $s \in [1,\infty)$ and further
\begin{align}\label{29-3-16ct1}
R(\widehat{q},\widehat{a}) \le \liminfn R(q_n,a_n).
\end{align}
Using Lemma \ref{weakly conv.} and the identities \eqref{27-9-19ct1}, we get that
\begin{align}\label{29-3-16ct2*}
J_\delta(\widehat{q},\widehat{a}) = \limn J^{h_n}_{\delta}\left(q_n,a_n\right).
\end{align}
Furthermore, since
$$\limn \|r^\epsilon_n - r^\epsilon\|_{L^1(\Omega)}
= \limn \|b^\epsilon_n - b^\epsilon\|_{L^1(\Omega)} = 0,$$
we also have
\begin{align}\label{29-3-16ct2}
J_\delta \left( r^\epsilon, b^\epsilon \right) = \limn J^{h_n}_{\delta}\left(r^\epsilon_n, b^\epsilon_n\right).
\end{align}
On the other hand, by the definition of $(r^\epsilon,b^\epsilon)$, we get 
$$|r^\epsilon -q| \le |q^\epsilon - q| \en \mbox{and} \en |b^\epsilon -a| \le |a^\epsilon - a|$$
 a.e. in $\Omega$. Integrating the above inequalities over the domain $\Omega$, it gives
$$\left\|r^\epsilon - q\right\|_{L^1(\Omega)} + \left\|b^\epsilon - a\right\|_{L^1(\Omega)} \le \left\|q^\epsilon - q\right\|_{L^1(\Omega)} + \left\|a^\epsilon - a\right\|_{L^1(\Omega)} \le C\epsilon$$
together with the limit
\begin{align}\label{29-3-16ct2**}
J_\delta(q,a) = \lim_{\epsilon\to 0}J_\delta \left( r^\epsilon, b^\epsilon \right).
\end{align}
We mention that
\begin{align*}
\left\|r^\epsilon - q\right\|_{L^2(\Omega)} + \left\|b^\epsilon - a\right\|_{L^2(\Omega)}
&\le 2 \max(\overline{q},\overline{a})^{1/2} \left(\left\|r^\epsilon - q\right\|_{L^1(\Omega)} + \left\|b^\epsilon - a\right\|_{L^1(\Omega)} \right)^{1/2}\\
&\le C \epsilon^{1/2}
\end{align*}
and then
\begin{align*}
\left\|r^\epsilon \right\|^2_{L^2(\Omega)} + \left\|b^\epsilon\right\|^2_{L^2(\Omega)}
\le C\epsilon + \left\|q \right\|^2_{L^2(\Omega)} + \left\|a\right\|^2_{L^2(\Omega)}.
\end{align*}
Combining this with \eqref{6-4-16ct5*}, it gives
\begin{align}\label{12-9-19ct2}
R(r^\epsilon,b^\epsilon) \le C\epsilon + R(q,a).
\end{align}
Furthermore, with the aid of \eqref{16-11-18ct1}, we get
\begin{align}\label{12-9-19ct3}
R(r^\epsilon,b^\epsilon) = \limn R(r^\epsilon_n,b^\epsilon_n).
\end{align}
Therefore, we obtain from by \eqref{29-3-16ct1}, \eqref{29-3-16ct2*}, \eqref{29-8-19ct1}, \eqref{29-3-16ct2}, \eqref{12-9-19ct3} and \eqref{12-9-19ct2} that
\begin{align*}
\Upsilon_{\delta,\rho}(\widehat{q},\widehat{a}) 
&= J_\delta(\widehat{q},\widehat{a}) + \rho R(\widehat{q},\widehat{a}) \notag\\
&\le \limn J^{h_n}_{\delta}\left(q_n,a_n\right) + \liminfn \rho R(q_n,a_n) \notag\\
&=\liminfn \left(J^{h_n}_{\delta}\left(q_n,a_n\right) + \rho R(q_n,a_n)\right)\notag\\
&\le \liminfn \left(  J^{h_n}_{\delta} \left(r^\epsilon_n,b^\epsilon_n\right)  + \rho R(r^\epsilon_n,b^\epsilon_n)\right),  \notag\\
&= J_\delta \left( r^\epsilon,b^\epsilon \right) + \rho R (r^\epsilon, b^\epsilon),\notag\\
&\le J_\delta \left( r^\epsilon,b^\epsilon \right) + \rho R(q,a) + C\epsilon\rho.
\end{align*}
Sending $\epsilon \to 0$, by \eqref{29-3-16ct2**}, we arrive at
\begin{align*}
\Upsilon_{\delta,\rho}(\widehat{q},\widehat{a}) 
\le J_\delta \left( q,a \right) + \rho R(q,a).
\end{align*}
Since $(q,a)$ is arbitrarily taken  in the admissible set $Q_{ad} \times A_{ad}$, the last relation shows that $(\widehat{q},\widehat{a})$ is a solution to
$\left( \mathcal{P}_{\rho,\delta} \right)$.

Now, denoting
$$\big(\widehat{r}^\epsilon, \widehat{b}^\epsilon \big) := \big(\widehat{r}^\epsilon (\widehat{q}), \widehat{b}^\epsilon (\widehat{a}) \big)
\quad \mbox{and} \quad
\big(\widehat{r}^\epsilon_n, \widehat{b}^\epsilon_n \big)
:= \big(I^{h_n}_1\widehat{r}^\epsilon, I^{h_n}_1\widehat{b}^\epsilon \big),$$
we have 
\begin{align*}
\rho\limsupn R(q_n,a_n)
&= \limn J^{h_n}_{\delta}\left(q_n,a_n\right) + \rho\limsupn R(q_n,a_n) - J_\delta(\widehat{q},\widehat{a})\\
&= \limsupn \left( J^{h_n}_{\delta}\left(q_n,a_n\right) + \rho R(q_n,a_n) \right)  - J_\delta(\widehat{q},\widehat{a}) \\
&\le \limsupn \left( J^{h_n}_{\delta} \big(\widehat{r}^\epsilon_n,\widehat{b}^\epsilon_n\big) + \rho R\big(\widehat{r}^\epsilon_n, \widehat{b}^\epsilon_n\big)\right)   - J_\delta(\widehat{q},\widehat{a})\\
&= \limn J^{h_n}_{\delta} \big(\widehat{r}^\epsilon_n,\widehat{b}^\epsilon_n\big) + \rho \limn R\big(\widehat{r}^\epsilon_n, \widehat{b}^\epsilon_n\big)   - J_\delta(\widehat{q},\widehat{a})\\
&= J_\delta\big(\widehat{r}^\epsilon,\widehat{b}^\epsilon\big) + \rho R \big(\widehat{r}^\epsilon,\widehat{b}^\epsilon\big)  - J_\delta(\widehat{q},\widehat{a})\\
&\le J_\delta\big(\widehat{r}^\epsilon,\widehat{b}^\epsilon\big) + \rho R (\widehat{q},\widehat{a}) + C\epsilon\rho - J_\delta(\widehat{q},\widehat{a}).
\end{align*}
Sending $\epsilon \to 0$, we get
\begin{align*}
\rho\limsupn R(q_n,a_n)
\le J_\delta(\widehat{q},\widehat{a}) + \rho R (\widehat{q},\widehat{a})  - J_\delta(\widehat{q},\widehat{a}) = \rho R (\widehat{q},\widehat{a}).
\end{align*}
This together with \eqref{29-3-16ct1} infers 
\begin{align*}
R(\widehat{q},\widehat{a}) \le \liminfn R(q_n,a_n) \le \limsupn R(q_n,a_n) \le R (\widehat{q},\widehat{a})
\end{align*}
and thus
\begin{align*}
\limn \left( \int_\Omega |\nabla q_n| + \int_\Omega |\nabla a_n| \right) = \int_\Omega |\nabla \widehat{q}| + \int_\Omega |\nabla \widehat{a}|.
\end{align*}
Utilizing Lemma \ref{bv1} again, we have
\begin{align*}
\int_\Omega |\nabla \widehat{q}| 
& \le \liminfn \int_\Omega |\nabla q_n|\\
&= \limn \left( \int_\Omega |\nabla q_n| + \int_\Omega |\nabla a_n| \right) - \liminfn \int_\Omega |\nabla a_n|\\
&= \int_\Omega |\nabla \widehat{q}| + \int_\Omega |\nabla \widehat{a}| - \liminfn \int_\Omega |\nabla a_n|
\end{align*}
and arrive at 
\begin{align*}
\liminfn \int_\Omega |\nabla a_n| \le \int_\Omega |\nabla \widehat{a}| \le \liminfn \int_\Omega |\nabla a_n|.
\end{align*}
This leads to the identity \eqref{28-8-19ct2}. Finally, since $\left(q_n,a_n\right)$ converges to $(\widehat{q},\widehat{a})$ in the $L^1(\Omega)$-norm and \eqref{28-8-19ct2}, we conclude 
that $(q_n, a_n)$
weakly converges to $(\widehat{q},\widehat{a})$ in $BV(\Omega) \times BV(\Omega)$ (see \cite{attouch},
Proposition 10.1.2, p. 374). Therefore, \eqref{29-8-19ct4} follows. The theorem is proved.
\end{proof}

We now introduce the notion of the unique $TV-L^2$-minimizing solution of the identification problem.

\begin{lemma}\label{10-9-19ct1}
The problem 
$$
\min_{\left\{(q,a)\in Q_{ad}\times A_{ad} ~\big|~ N_{j^\dag}(q,a) = M_{g^\dag}(q,a) \right\}} R(q,a) \eqno\left(\mathcal{IP}\right)
$$
admits a solution, which is called the $TV-L^2$-minimizing solution of the identification problem.
\end{lemma}

\begin{proof}
The assertion follows from standard arguments, it is therefore ignored here.
\end{proof}

\begin{lemma}\label{23-6-16ct5}
For any fixed $(q,a)\in Q_{ad}\times A_{ad}$ an element $(\widehat{q}^h,\widehat{a}^h) \in Q^h_{ad} \times A^h_{ad}$ exists such that
\begin{align}\label{22-7-16ct5}
\big\|\widehat{q}^h-q\big\|_{L^1(\Omega)} + \big\|\widehat{a}^h-a\big\|_{L^1(\Omega)} \le C h|\log h|
\end{align}
and
\begin{align}\label{22-7-16ct6}
\lim_{h\to 0}R (\widehat{q}^h,\widehat{a}^h) = R (q,a).
\end{align}
\end{lemma}

\begin{proof}
The existence of the pair $(\widehat{q}^h,\widehat{a}^h) \in Q^h_{ad} \times A^h_{ad}$ satisfying the inequality \eqref{22-7-16ct5} follows from Lemma 4.6 of \cite{HKQ18}, where
\begin{align}\label{28-1-19ct1}
\lim_{h\to0}\int_\Omega|\nabla \widehat{q}^h| = \int_\Omega|\nabla q| \quad\mbox{and}\quad \lim_{h\to0}\int_\Omega|\nabla \widehat{a}^h| = \int_\Omega|\nabla a|.
\end{align}
Since $\lim_{h\to0} h|\log h| =0$, the identity \eqref{22-7-16ct6} is now implied by \eqref{28-1-19ct1}. The proof completes.
\end{proof}

For any $(q,a)\in Q_{ad}\times A_{ad}$ let $(\widehat{q}^h,\widehat{a}^h) \in Q^h_{ad} \times A^h_{ad}$ be arbitrarily generated   from $(q,a)$. We have the limit
\begin{align*}
\chi^h_{j_\delta,g_\delta}(q,a) := \left\| N^h_{j_\delta}(\widehat{q}^h,\widehat{a}^h) - N_{j_\delta}(q,a) \right\|_{H^1(\Omega)} + \left\| M^h_{g_\delta}(\widehat{q}^h,\widehat{a}^h) - M_{g_\delta}(q,a) \right\|_{H^1(\Omega)} \rightarrow 0 \mbox{~as~} h\to 0
\end{align*}
and the estimate
\begin{equation}\label{eq:r}
\chi^h_{j_\delta,g_\delta}(q,a) \le C_r \big( h|\log h|\big)^r \quad\mbox{with}\quad
\begin{cases}
r<1/2&  \mbox{if} \quad d=2 \quad\mbox{and}\\
r=1/3&  \mbox{if} \quad d=3
\end{cases}
\end{equation}
in case $N_{j_\delta}(q,a),~ M_{g_\delta}(q,a) \in H^2(\Omega)$ (see \cite[Lemma 4.8]{HKQ18}).

\begin{theorem}\label{convergence1}
Let $\left(h_n\right)$,  $\left(\delta_n\right)$
and $\left(\rho_n\right)$ be any positive sequences such that
\begin{align}\label{29-6-16ct1}
\rho_n\rightarrow 0, ~\frac{\delta_n}{\sqrt{\rho_n}}
\rightarrow 0 \mbox{~and~} \frac{\chi^{h_n}_{j^\dag,g^\dag}(q,a)}{\sqrt{\rho_n}}
\rightarrow 0 \mbox{~as~} n\to\infty,
\end{align}
where $(q,a)$ is any solution of $N_{j^\dag}(q,a) = M_{g^\dag}(q,a)$. Moreover, assume that $\big(j_{\delta_n}, g_{\delta_n}\big) \subset H^{-1/2}(\Gamma) \times  H^{1/2}(\Gamma)$ is a sequence satisfying
$$\big\|j_{\delta_n} - j^\dag\big\|_{H^{-1/2}(\Gamma)} + \big\|g_{\delta_n} - g^\dag\big\|_{H^{1/2}(\Gamma)} \le \delta_n$$
and that $(q_n,a_n) := \big(q_{\rho_n, \delta_n}^{h_n},a_{\rho_n, \delta_n}^{h_n}\big)$ is an arbitrary minimizer of $\big( \mathcal{P}_{\rho_n,\delta_n}^{h_n} \big)$ for each $n\in\N$. Then,

(i) There exist a subsequence of $(q_n, a_n)$ denoted by the same symbol and a solution $(q^\dag,a^\dag)$ to $\left(\mathcal{IP}\right)$ such that $(q_n,a_n)$ converges to $(q^\dag,a^\dag)$ in the $L^s(\Omega)$-norm for all $s \in [1,\infty)$ and
\begin{align}
&\limn  \int_\Omega |\nabla q_n|  = \int_\Omega |\nabla q^\dag|  \quad\mbox{and}\quad \limn  \int_\Omega |\nabla a_n|  = \int_\Omega |\nabla a^\dag|,\label{28-8-19ct2*}\\
&\limn D_{TV}^{\ell}(q_n,q^\dag) = \limn D_{TV}^{\kappa}(a_n,a^\dag) = 0 \label{29-8-19ct4*}
\end{align}
for all
$(\ell,\kappa) \in \partial\left(\int_\Omega|\nabla(\cdot)|\right)(q^\dag) \times \partial\left(\int_\Omega|\nabla(\cdot)|\right)(a^\dag)$.

(ii) The sequences $\big( N^{h_n}_{j_{\delta_n}}(q_n,a_n)\big) $ and $\big( M^{h_n}_{g_{\delta_n}}(q_n,a_n)\big) $ converge in the $H^1(\Omega)$-norm to the unique weak solution $u(q^\dag,a^\dag)$ of the boundary value problem \eqref{17-5-16ct1} -- \eqref{17-5-16ct1***}.
\end{theorem}

\begin{proof}
We have from the equation $N_{j^\dag}(q,a) = M_{g^\dag}(q,a)$ and the optimality of $(q_n,a_n)$ that
\begin{align}\label{eq:opt}
J^{h_n}_{\delta_n} \left(q_n,a_n\right) + \rho_n R\left(q_n,a_n\right)
&\le J^{h_n}_{\delta_n} \big({\widehat{q}}^{h_n} ,{\widehat{a}}^{h_n} \big) + \rho_n R\big({\widehat{q}}^{h_n} ,{\widehat{a}}^{h_n} \big),
\end{align}
where $\big({\widehat{q}}^{h_n} ,{\widehat{a}}^{h_n} \big)$ is generated from $(q,a)$ according to Lemma \ref{23-6-16ct5}, and
\begin{align*}
J^{h_n}_{\delta_n} \big({\widehat{q}}^{h_n} ,{\widehat{a}}^{h_n} \big)
&\le C \left\| N^{h_n}_{j_{\delta_n}} ({\widehat{q}}^{h_n} , {\widehat{a}}^{h_n} ) -
M^{h_n}_{g_{\delta_n}} ({\widehat{q}}^{h_n} , {\widehat{a}}^{h_n} ) \right\|^2_{H^1(\Omega)}\\
&\le C \Bigl( \left\|N^{h_n}_{j_{\delta_n}} ({\widehat{q}}^{h_n} , {\widehat{a}}^{h_n} ) - N^{h_n}_{j^\dag} ({\widehat{q}}^{h_n} , {\widehat{a}}^{h_n} ) \right\|^2_{H^1(\Omega)}
+ \left\| M^{h_n}_{g^\dag} ({\widehat{q}}^{h_n} , {\widehat{a}}^{h_n} ) - M^{h_n}_{g_{\delta_n}} ({\widehat{q}}^{h_n} , {\widehat{a}}^{h_n} ) \right\|^2_{H^1(\Omega)}\\
&~\qquad + \left\|N^{h_n}_{j^\dag} ({\widehat{q}}^{h_n} , {\widehat{a}}^{h_n} ) - N_{j^\dag} (q,a)  \right\|^2_{H^1(\Omega)}
+  \left\|  M^{h_n}_{g^\dag} ({\widehat{q}}^{h_n} , {\widehat{a}}^{h_n} ) - M_{g^\dag} (q,a)  \right\|^2_{H^1(\Omega)}\Bigr)\\
&\le C \left( \left\| j_{\delta_n} - j^\dag \right\|^2_{H^{-1/2}(\Gamma)} + \left\| g_{\delta_n} - g^\dag \right\|^2_{H^{1/2}(\Gamma)}\right) + C \chi^{h_n}_{j^\dag,g^\dag}(q,a)^2\\
&\le C\left( \delta^2_n + \chi^{h_n}_{j^\dag,g^\dag}(q,a)^2\right).
\end{align*}
Therefore it follows from \eqref{eq:opt}, \eqref{29-6-16ct1} and \eqref{22-7-16ct6} that
\begin{align}\label{30-3-16ct2}
\limn J^{h_n}_{\delta_n} \left( q_n,a_n\right) =0
\end{align}
and
\begin{align}\label{30-3-16ct3}
\limsupn R(q_n,a_n)
\le \limsupn R\big({\widehat{q}}^{h_n} ,{\widehat{a}}^{h_n} \big)
= R(q,a).
\end{align}
With the aid of Lemma \ref{bv1}, a subsequence of $(q_n,a_n)$ not relabeled and an element $(q^\dag,a^\dag)\in Q_{ad} \times A_{ad}$ exist such that $\left(q_n,a_n\right)$ converges to $(q^\dag,a^\dag)$ in the $L^s(\Omega)$-norm for all $s \in [1,\infty)$ and
\begin{align}\label{30-3-16ct4}
R(q^\dag,a^\dag) \le \liminfn R(q_n,a_n).
\end{align}
Thus, due to Lemma \ref{weakly conv.}, we obtain that $\big( N^{h_n}_{j_{\delta_n}}(q_n,a_n),~ M^{h_n}_{g_{\delta_n}}(q_n,a_n)\big) $ converges to $\big(N_{j^\dag}\big(q^\dag,a^\dag),~ M_{g^\dag}\big(q^\dag,a^\dag\big)\big)$ in the $H^1(\Omega)\times H^1(\Omega)$-norm. This yields the equation
\begin{align*}
\big\|N_{j^\dag}\big(q^\dag,a^\dag\big) - M_{g^\dag}\big(q^\dag,a^\dag\big)\big\|_{H^1(\Omega)}
&=\limn \big\|N^{h_n}_{j_{\delta_n}}(q_n,a_n) - M^{h_n}_{g_{\delta_n}}(q_n,a_n)\big\|_{H^1(\Omega)}\\
&\le C \limn \sqrt{ J^{h_n}_{\delta_n} (q_n,a_n)} =0,
\end{align*}
by \eqref{30-3-16ct2}. Thus, $(q^\dag,a^\dag)$ belongs to the set $\left\{(q,a)\in Q_{ad}\times A_{ad} ~\big|~ N_{j^\dag}(q,a) = M_{g^\dag}(q,a) \right\}$.

It follows from \eqref{30-3-16ct3} -- \eqref{30-3-16ct4} that
\begin{align}\label{8-1-20ct1}
R \big(q^\dag,a^\dag) \le \liminfn R (q_n,a_n) \le \limsupn  R (q_n,a_n) \le  R (q,a)
\end{align}
and $(q^\dag,a^\dag)$ is thus a solution to $\left(\mathcal{IP}\right)$. Further, replacing $(q,a)$ in \eqref{8-1-20ct1} by $(q^\dag,a^\dag)$, we obtain
$$R \big(q^\dag,a^\dag) = \limn R (q_n,a_n).$$
Therefore, using the arguments included in the proof of Theorem \ref{29-8-19ct5}, we arrive at \eqref{28-8-19ct2*} and \eqref{29-8-19ct4*}, which finishes the proof.
\end{proof}

\section{Differential and projected gradient algorithm}\label{29-8-19ct3}

We start the section with presenting the differentials of the discrete coefficient-to-solution operators and of the associated cost functional.

\begin{lemma}\label{qh-diff} The discrete operators $N^h_{j_\delta}$ and $M^h_{g_\delta}$ are infinitely Fr\'echet differentiable. For $(q,a) \in Q\times A$ and $\big((\eta_q^1,\eta_a^1), \ldots, (\eta_q^m,\eta_a^m) \big) \in {L^\infty(\Omega)}^{2m}$, the $m$-th order differentials 
$${D^h_N}^{(m)} := {N_{j_\delta}^h}^{(m)}(q,a)\big((\eta_q^1,\eta_a^1), \ldots, (\eta_q^m,\eta_a^m) \big) \in V^h_1$$
and
$${D^h_M}^{(m)} := {M^h_{g_\delta}}^{(m)}(q,a)\big((\eta_q^1,\eta_a^1), \ldots, (\eta_q^m,\eta_a^m) \big) \in V^h_{1,0}$$ are the unique solutions to the variational equations
\begin{align*}
&\int_\Omega q\nabla {D^h_N}^{(m)} \cdot \nabla \phi^h dx  
+ \int_\Omega a {D^h_N}^{(m)}  \phi^h dx + \int_{\partial\Omega} \sigma {D^h_N}^{(m)} \phi^h ds \notag\\
&~\quad = -\sum_{i=1}^m \int_\Omega \eta_q^i \nabla {N^h_{j_\delta}}^{(m-1)}(q,a) \overline{\eta}_q^i \cdot \nabla\phi^h dx -\sum_{i=1}^m \int_\Omega \eta_a^i {N^h_{j_\delta}}^{(m-1)}(q,a) \overline{\eta}_a^i\phi^h dx, \quad \forall \phi^h\in V^h_1
\end{align*}
and
\begin{align*}
&\int_\Omega q\nabla {D^h_M}^{(m)} \cdot \nabla \phi^h dx  
+ \int_\Omega a {D^h_M}^{(m)}  \phi^h dx + \int_{\partial\Omega} \sigma {D^h_M}^{(m)} \phi^h ds \notag\\
&~\quad = -\sum_{i=1}^m \int_\Omega \eta_q^i \nabla {M^h_{g_\delta}}^{(m-1)}(q,a) \overline{\eta}_q^i \cdot \nabla\phi^h dx -\sum_{i=1}^m \int_\Omega \eta_a^i {M^h_{g_\delta}}^{(m-1)}(q,a) \overline{\eta}_a^i\phi^h dx, \quad \forall \phi^h\in V^h_{1,0}
\end{align*}
with $\overline{\eta}_q^i := (\eta_q^1, \ldots, \eta_q^{i-1},\eta_q^{i+1}, \ldots, \eta_q^m) \in {L^\infty(\Omega)}^{m-1}$ and $\overline{\eta}_a^i := (\eta_a^1, \ldots, \eta_a^{i-1},\eta_a^{i+1}, \ldots, \eta_a^m) \in {L^\infty(\Omega)}^{m-1}$, respectively. Furthermore, 
\begin{align*}
\max\left(\left\|{D^h_N}^{(m)}\right\|_{H^1(\Omega)}, \left\|{D^h_M}^{(m)}\right\|_{H^1(\Omega)} \right) \le C \prod_{i=1}^m \left( \|\eta_q^i\|_{L^\infty(\Omega)} + \|\eta_a^i\|_{L^\infty(\Omega)} \right).
\end{align*}
\end{lemma}

\begin{proof}
The proof is based on standard arguments, is therefore omitted here.
\end{proof}

Below we present the gradient of the cost functional. For $(q,a)\in Q^h_{ad} \times A^h_{ad}$ and $(\eta_q,\eta_a) \in V^h_1 \times V^h_1$ we get that
\begin{align*}
{J^h_\delta}'(q,a)(\eta_q,\eta_a) = \frac{\partial J^h_\delta(q,a)}{\partial q}\eta_q + \frac{\partial J^h_\delta(q,a)}{\partial a}\eta_a,
\end{align*}
where
\begin{align*}
\frac{1}{2}\frac{\partial J^h_\delta(q,a)}{\partial q}\eta_q 
&= \frac{1}{2}\int_\Omega \eta_q \left|\nabla\big(N^h_{j_\delta}(q,a) - M^h_{g_\delta}(q,a)\big)\right|^2 dx \\
&~\quad + \int_\Omega q \nabla\big(N^h_{j_\delta}(q,a) - M^h_{g_\delta}(q,a)\big) \cdot \nabla\big({N^h_{j_\delta}}'(q,a)(\eta_q,0) - {M^h_{g_\delta}}'(q,a)(\eta_q,0)\big)dx\\
&~\quad + \int_\Omega a\big(N^h_{j_\delta}(q,a) - M^h_{g_\delta}(q,a)\big) \big({N^h_{j_\delta}}'(q,a)(\eta_q,0) - {M^h_{g_\delta}}'(q,a)(\eta_q,0)\big) dx
\\
&~\quad + \int_{\partial\Omega} \sigma \big(N^h_{j_\delta}(q,a) - M^h_{g_\delta}(q,a)\big) \big({N^h_{j_\delta}}'(q,a)(\eta_q,0) - {M^h_{g_\delta}}'(q,a)(\eta_q,0)\big) ds
\end{align*}
and 
\begin{align*}
\frac{1}{2}\frac{\partial J^h_\delta(q,a)}{\partial a}\eta_a 
&= \int_\Omega q \nabla\big(N^h_{j_\delta}(q,a) - M^h_{g_\delta}(q,a)\big) \cdot \nabla\big({N^h_{j_\delta}}'(q,a)(0,\eta_a) - {M^h_{g_\delta}}'(q,a)(0,\eta_a)\big)dx\\
&~\quad + \frac{1}{2}\int_\Omega \eta_a \big(N^h_{j_\delta}(q,a) - M^h_{g_\delta}(q,a)\big)^2 dx\\
&~\quad + \int_\Omega a\big(N^h_{j_\delta}(q,a) - M^h_{g_\delta}(q,a)\big) \big({N^h_{j_\delta}}'(q,a)(0,\eta_a) - {M^h_{g_\delta}}'(q,a)(0,\eta_a)\big) dx \\
&~\quad + \int_{\partial\Omega} \sigma \big(N^h_{j_\delta}(q,a) - M^h_{g_\delta}(q,a)\big) \big({N^h_{j_\delta}}'(q,a)(0,\eta_a) - {M^h_{g_\delta}}'(q,a)(0,\eta_a)\big) ds.
\end{align*}
Thus,
\begin{align*}
\frac{1}{2}{J^h_\delta}'(q,a)(\eta_q,\eta_a) 
&= \frac{1}{2}\int_\Omega \eta_q \left|\nabla\big(N^h_{j_\delta}(q,a) - M^h_{g_\delta}(q,a)\big)\right|^2 dx + \frac{1}{2}\int_\Omega \eta_a \big(N^h_{j_\delta}(q,a) - M^h_{g_\delta}(q,a)\big)^2 dx\\
&~\quad + \int_\Omega q \nabla\big(N^h_{j_\delta}(q,a) - M^h_{g_\delta}(q,a)\big) \cdot \nabla\big({N^h_{j_\delta}}'(q,a)(\eta_q,\eta_a) - {M^h_{g_\delta}}'(q,a)(\eta_q,\eta_a)\big)dx\\
&~\quad + \int_\Omega a\big(N^h_{j_\delta}(q,a) - M^h_{g_\delta}(q,a)\big) \big({N^h_{j_\delta}}'(q,a)(\eta_q,\eta_a) - {M^h_{g_\delta}}'(q,a)(\eta_q,\eta_a)\big) dx
\\
&~\quad + \int_{\partial\Omega} \sigma \big(N^h_{j_\delta}(q,a) - M^h_{g_\delta}(q,a)\big) \big({N^h_{j_\delta}}'(q,a)(\eta_q,\eta_a) - {M^h_{g_\delta}}'(q,a)(\eta_q,\eta_a)\big) ds.
\end{align*}
Denoting by $\sum$ the last three terms in the above sum, we have 
\begin{align*}
\sum &:=  \int_\Omega q \nabla {N^h_{j_\delta}}'(q,a)(\eta_q,\eta_a)\cdot \nabla\big(N^h_{j_\delta}(q,a) - M^h_{g_\delta}(q,a)\big)   dx \\
&~\quad + \int_\Omega a {N^h_{j_\delta}}'(q,a)(\eta_q,\eta_a)  \big(N^h_{j_\delta}(q,a) - M^h_{g_\delta}(q,a)\big) dx
\\
&~\quad + \int_{\partial\Omega} \sigma {N^h_{j_\delta}}'(q,a)(\eta_q,\eta_a)\big(N^h_{j_\delta}(q,a) - M^h_{g_\delta}(q,a)\big)  ds\\
&~\quad + \int_\Omega q \nabla M^h_{g_\delta}(q,a) \cdot \nabla {M^h_{g_\delta}}'(q,a)(\eta_q,\eta_a)dx + \int_\Omega a M^h_{g_\delta}(q,a) {M^h_{g_\delta}}'(q,a)(\eta_q,\eta_a) dx
\\
&~\quad + \int_{\partial\Omega} \sigma M^h_{g_\delta}(q,a)\big)  {M^h_{g_\delta}}'(q,a)(\eta_q,\eta_a) ds \\
&~\quad - \int_\Omega q \nabla N^h_{j_\delta}(q,a) \cdot \nabla {M^h_{g_\delta}}'(q,a)(\eta_q,\eta_a)dx - \int_\Omega a N^h_{j_\delta}(q,a) {M^h_{g_\delta}}'(q,a)(\eta_q,\eta_a) dx
\\
&~\quad - \int_{\partial\Omega} \sigma N^h_{j_\delta}(q,a)\big)  {M^h_{g_\delta}}'(q,a)(\eta_q,\eta_a) ds.
\end{align*}
With the aid of Lemma \ref{qh-diff} together with \eqref{10/4:ct1} and \eqref{10/4:ct1*} we get
\begin{align*}
\sum 
& = - \int_\Omega \eta_q \nabla {N^h_{j_\delta}}(q,a)\cdot \nabla\big(N^h_{j_\delta}(q,a) - M^h_{g_\delta}(q,a)\big)   dx - \int_\Omega \eta_a {N^h_{j_\delta}}(q,a) \big(N^h_{j_\delta}(q,a) - M^h_{g_\delta}(q,a)\big) dx \\
&~\quad + \<f,{M^h_{g_\delta}}'(q,a)(\eta_q,\eta_a)\>_{\big(H^{-1}(\Omega), H^1(\Omega) \big)}  + \<j_0, {M^h_{g_\delta}}'(q,a)(\eta_q,\eta_a)\>_{\big( H^{-1/2}(\partial\Omega\setminus\Gamma), H^{1/2}(\partial\Omega\setminus\Gamma) \big)} \\
&~\quad - \<f,{M^h_{g_\delta}}'(q,a)(\eta_q,\eta_a)\>_{\big(H^{-1}(\Omega), H^1(\Omega) \big)} - \<j_\delta, {M^h_{g_\delta}}'(q,a)(\eta_q,\eta_a)\>_{\big( H^{-1/2}(\Gamma), H^{1/2}(\Gamma) \big)} \\
&~\quad - \<j_0, {M^h_{g_\delta}}'(q,a)(\eta_q,\eta_a)\>_{\big( H^{-1/2}(\partial\Omega\setminus\Gamma), H^{1/2}(\partial\Omega\setminus\Gamma) \big)} \\
&= - \int_\Omega \eta_q \nabla {N^h_{j_\delta}}(q,a)\cdot \nabla\big(N^h_{j_\delta}(q,a) - M^h_{g_\delta}(q,a)\big)   dx - \int_\Omega \eta_a {N^h_{j_\delta}}(q,a) \big(N^h_{j_\delta}(q,a) - M^h_{g_\delta}(q,a)\big) dx,
\end{align*}
due to the fact ${M^h_{g_\delta}}'(q,a)(\eta_q,\eta_a) \in V^h_{1,0}$. Consequently, we obtain that
\begin{align*}
{J^h_\delta}'(q,a)(\eta_q,\eta_a) 
&= \int_\Omega \eta_q \left|\nabla\big(N^h_{j_\delta}(q,a) - M^h_{g_\delta}(q,a)\big)\right|^2 dx +\int_\Omega \eta_a \big(N^h_{j_\delta}(q,a) - M^h_{g_\delta}(q,a)\big)^2 dx\\
&~\quad - 2\int_\Omega \eta_q \nabla {N^h_{j_\delta}}(q,a)\cdot \nabla\big(N^h_{j_\delta}(q,a) - M^h_{g_\delta}(q,a)\big)   dx \\
&~\quad - 2\int_\Omega \eta_a {N^h_{j_\delta}}(q,a) \big(N^h_{j_\delta}(q,a) - M^h_{g_\delta}(q,a)\big) dx.
\end{align*}
Therefore, we arrive at the following result.

\begin{lemma}\label{2-9-19ct1}
The differential of the functional $J^h_\delta$ at $(q,a)\in Q^h_{ad}\times A^h_{ad}$ in the direction $(\eta_q,\eta_a) \in V^h_1 \times V^h_1$ given by
\begin{align}\label{2-9-19ct3}
{J^h_\delta}'(q,a)(\eta_q,\eta_a) =  \int_\Omega \eta_q \left( \left| \nabla {M^h_{g_\delta}}(q,a) \right|^2 - \left| \nabla {N^h_{j_\delta}}(q,a) \right|^2 \right)dx +  \int_\Omega \eta_a \left( \left| {M^h_{g_\delta}}(q,a) \right|^2 - \left| {N^h_{j_\delta}}(q,a) \right|^2 \right)dx.
\end{align}
\end{lemma}

The mainly computational challenge of the total variation regularization method is non-differentiable of the $BV$-semi-norm. To overcome this difficulty, we replace the total variation by a differentiable approximation
\begin{align*}
\int_\Omega |\nabla q| \approx \int_\Omega \sqrt{|\nabla q|^2 + \epsilon^h} dx \quad\mbox{and}\quad \int_\Omega |\nabla a| \approx \int_\Omega \sqrt{|\nabla a|^2 + \epsilon^h} dx,
\end{align*}
where $\epsilon^h$ is a positive function of the mesh size $h$ satisfying $\lim_{h\to 0} \epsilon^h =0$. Thus, the regularization term is approximated by
\begin{align*}
R(q,a) \approx R^\epsilon(q,a) := \int_\Omega \sqrt{|\nabla q|^2 + \epsilon^h} dx + \int_\Omega \sqrt{|\nabla a|^2 + \epsilon^h} dx + \frac{1}{2}\|q\|^2_{L^2(\Omega)} + \frac{1}{2} \|a\|^2_{L^2(\Omega)}.
\end{align*}
For all $(\eta_q,\eta_a) \in V^h_1 \times V^h_1$ we get that
\begin{align}\label{2-9-19ct4}
{R^\epsilon}'(q,a)(\eta_q,\eta_a) 
= \int_\Omega \frac{\nabla \eta_q \cdot \nabla q}{\sqrt{|\nabla q|^2 + \epsilon^h}} dx + \int_\Omega \frac{\nabla \eta_a \cdot \nabla a}{\sqrt{|\nabla a|^2 + \epsilon^h}} dx + \int_\Omega \eta_q q dx + \int_\Omega \eta_a a dx.
\end{align}
The discrete cost functional $\Upsilon^h_{\delta,\rho} (q,a)$ of the problem $\big(\mathcal{P}^h_{\delta,\rho}\big)$ is then approximated by
\begin{align*}
\Upsilon^{h,\epsilon}_{\delta,\rho} (q,a) := J^h_\delta(q,a) + \rho R^\epsilon(q,a).
\end{align*}

\begin{lemma}\label{2-9-19ct2}
The differential of the approximated cost functional $\Upsilon^{h,\epsilon}_{\delta,\rho}$ at $(q,a)\in Q^h_{ad}\times A^h_{ad}$ in the direction $(\eta_q,\eta_a) \in V^h_1 \times V^h_1$ fulfilled the identity
\begin{align}\label{2-6-20ct2}
{\Upsilon^{h,\epsilon}_{\delta,\rho}}' (q,a)(\eta_q,\eta_a) 
&= \int_\Omega \eta_q \left( \left| \nabla {M^h_{g_\delta}}(q,a) \right|^2 - \left| \nabla {N^h_{j_\delta}}(q,a) \right|^2 \right)dx +  \int_\Omega \eta_a \left( \left| {M^h_{g_\delta}}(q,a) \right|^2 - \left| {N^h_{j_\delta}}(q,a) \right|^2 \right)dx \notag\\
&~\quad + \rho\left( \int_\Omega \frac{\nabla \eta_q \cdot \nabla q}{\sqrt{|\nabla q|^2 + \epsilon^h}} dx + \int_\Omega \frac{\nabla \eta_a \cdot \nabla a}{\sqrt{|\nabla a|^2 + \epsilon^h}} dx + \int_\Omega \eta_q q dx + \int_\Omega \eta_a a dx\right).
\end{align}
\end{lemma}

\begin{proof}
The affirmation directly follows  from the definition of the functional  $\Upsilon^{h,\epsilon}_{\delta,\rho}$ and the identities \eqref{2-9-19ct3} -- \eqref{2-9-19ct4}.
\end{proof}

With the derivative ${\Upsilon^{h,\epsilon}_{\delta,\rho}}' (q,a)$ of the approximated cost functional $\Upsilon^{h,\epsilon}_{\delta,\rho}$ at $(q,a)\in Q^h_{ad}\times A^h_{ad}$ at hand, we now present a projected gradient method to reach a minimizer to $(\mathcal{P}^h_{\delta,\rho})$ (see  \cite{blank} and the references given there for detailed discussions on the method).

\begin{algorithm}[H] \label{31-8-18ct1}
    \SetKwInOut{Input}{Input}
    \SetKwInOut{Output}{Output}
     \vspace{0.3cm}
    \Input{Given an initial approximation $(q_0,a_0) \in Q^h_{ad} \times A^h_{ad}$, a smoothing parameter $\epsilon^h$, positive constants $\kappa_1,~ \kappa_2$, number of iteration $N$ and setting $n=0$.}
    \Output{An approximation of a solution to $(\mathcal{P}^h_{\delta,\rho})$.}
    \vspace{0.3cm}
    \While{$n\le N$}
      {
      1. Compute the gradient of the cost functional 
      $$\nabla \Upsilon^{h,\epsilon}_{\delta,\rho} (q_n,a_n) := \left( \nabla_q \Upsilon^{h,\epsilon}_{\delta,\rho} (q_n,a_n), \nabla_a \Upsilon^{h,\epsilon}_{\delta,\rho} (q_n,a_n)\right)$$
      
      2. Choose the maximum $\beta_n \in \{1,1/2,1/4,1/8, \ldots\}$ such that $\Upsilon^{h,\epsilon}_{\rho,\delta} (\widehat{q}_n,\widehat{a}_n) <  \Upsilon^{h,\epsilon}_{\rho,\delta} (q_n,a_n)$, where
      \begin{align*}
      \widehat{q}_n &= \max\left( \underline{q}, \min\left(q_n - \beta_n \nabla_q \Upsilon^{h,\epsilon}_{\delta,\rho} (q_n,a_n), \overline{q} \right) \right)\\
      \widehat{a}_n &= \max\left( \underline{a}, \min\left(a_n - \beta_n \nabla_a \Upsilon^{h,\epsilon}_{\delta,\rho} (q_n,a_n), \overline{a} \right) \right).
      \end{align*}
      
      3. Compute $ Tolerance := \big\| \nabla \Upsilon^h_{\rho,\delta}(\widehat{q}_n,\widehat{a}_n) \big\|_{L^2(\Omega)} -\kappa_1 -\kappa_2\big\| \nabla \Upsilon^h_{\rho,\delta}(q_0,a_0) \big\|_{L^2(\Omega)} $
      
         \eIf{$ Tolerance \le 0$}{
              stop
              }{
              set $n = n+1$ and update $(q_n,a_n) = (\widehat{q}_{n-1}, \widehat{a}_{n-1})$, then go back to Step 1.
              } 
      }
    \caption{Minimizing of $(\mathcal{P}^h_{\delta,\rho})$}
\end{algorithm}

In Step 1 of Algorithm \ref{31-8-18ct1} the gradient $\nabla \Upsilon^{h,\epsilon}_{\delta,\rho} = \left( \nabla_q \Upsilon^{h,\epsilon}_{\delta,\rho}, \nabla_a \Upsilon^{h,\epsilon}_{\delta,\rho} \right) \in V^h_1 \times V^h_1$ is given by
\begin{align}\label{2-6-20ct1}
\left( \nabla \Upsilon^{h,\epsilon}_{\delta,\rho} (q,a), (\eta_q,\eta_a) \right)_{V^h_1 \times V^h_1} = {\Upsilon^{h,\epsilon}_{\delta,\rho}}' (q,a)(\eta_q,\eta_a)
\end{align}
for all $(\eta_q,\eta_a) \in V^h_1 \times V^h_1$.  Let $\{\phi_1, \ldots,\phi_{E^h}\}$ be the basis of ${V}^h_1$ consisting hat functions, i.e. $\phi_i(N_j) =\delta_{ij}$ for all $1\le i,j\le E^h$, where $\delta_{ij}$ is the Kronecker symbol and $N_j$ is the $j^{th}$-node of the triangulation $\mathcal{T}^h$. Each function $\phi\in{V}^h_1$ can be then identified with a vector $(\phi_1,\ldots,\phi_{E^h})\in\mathbb{R}^{E^h}$ consisting of its nodal values, i.e.
$\phi=\sum_{j=1}^{E^h} \phi(N_j) \phi_j$. Now, denoting by $\nabla_q \Upsilon^{h,\epsilon}_{\delta,\rho}(q,a) = (\Upsilon_1^q(q,a),\ldots,\Upsilon_{E^h}^q(q,a))$ and taking $\eta_a \equiv 0$ in \eqref{2-6-20ct1} --    \eqref{2-6-20ct2}, we for each $j\in \{1,\ldots, E^h\}$ arrive at
\begin{align*}
\Upsilon^q_j(q,a)
= \int_\Omega \phi_j \left( \left| \nabla {M^h_{g_\delta}}(q,a) \right|^2 - \left| \nabla {N^h_{j_\delta}}(q,a) \right|^2 \right)dx
+ \rho\left( \int_\Omega \frac{\nabla \phi_j \cdot \nabla q}{\sqrt{|\nabla q|^2 + \epsilon^h}} dx + \int_\Omega \phi_j q dx \right).
\end{align*}
Likewise, with $\nabla_a \Upsilon^{h,\epsilon}_{\delta,\rho}(q,a) = (\Upsilon_1^a(q,a),\ldots,\Upsilon_{E^h}^a(q,a))$ it has
\begin{align*}
\Upsilon^a_j(q,a)
= \int_\Omega \phi_j \left( \left| {M^h_{g_\delta}}(q,a) \right|^2 - \left| {N^h_{j_\delta}}(q,a) \right|^2 \right)dx
+ \rho\left( \int_\Omega \frac{\nabla \phi_j \cdot \nabla a}{\sqrt{|\nabla a|^2 + \epsilon^h}} dx + \int_\Omega \phi_j a dx \right).
\end{align*}
In Step 2 the projected step size $\beta_n$ is chosen such that (cf.\ \cite[Chapter 2]{Hi09})
\begin{align}\label{2-6-20ct3}
\Upsilon^{h,\epsilon}_{\rho,\delta} (\widehat{q}_n,\widehat{a}_n) -  \Upsilon^{h,\epsilon}_{\rho,\delta} (q_n,a_n) 
\le -\frac{\beta}{\beta_n} \left( \| \widehat{q}_n -q_n\|^2_{V^h_1} + \| \widehat{a}_n -a_n\|^2_{V^h_1}\right)
\end{align}
for some $\beta \in (0,1)$.

\section{Numerical examples} \label{numerical examples}

Our numerical case study is the equation
\begin{align}
-\nabla \cdot \big(q^\dag \nabla u \big) + a^\dag u &= f \quad \mbox{~in} \quad \Omega,  \label{9-9-19ct1}\\
q^\dag \nabla u \cdot \vec{n} +\sigma u &= j^\dag \quad \mbox{on} \quad \Gamma,  \label{9-9-19ct2}\\
q^\dag \nabla u \cdot \vec{n} +\sigma u &= j_0 \quad \mbox{on} \quad \partial\Omega\setminus\Gamma,  \label{9-9-19ct3}\\
 u &= g^\dag \quad \mbox{on} \quad \Gamma  \label{9-9-19ct4}
\end{align}
with the domain $\Omega = \{ x = (x_1,x_2) \in \mathbb{R}^2 ~|~ -1 < x_1, x_2 < 1\}$ and the observation boundary $\Gamma := (-1,1) \times \{-1\} \cup \{-1\}\times [-1,1)$ (the bottom edge and the left edge).

The known functions are given as: the Robin coefficient $\sigma = 1$ on $\partial\Omega$, the Neumann data on $\partial\Omega\setminus\Gamma$
$$
j_0 := 4\chi_{\{1\}\times[-1,1]} - 3\chi_{[-1,1)\times\{1\}},
$$
and the source term
$$f := \chi_{D} - \chi_{\Omega \setminus D},$$ 
where $\chi_D$ is the characteristic function of the Lebesgue measurable set 
$$D := \left\{ (x_1, x_2) \in \Omega ~\big|~ |x_1| + |x_2| \le 1/2\right\}.$$
The sought diffusion and reaction coefficients $q^\dag$ and $a^\dag$ are respectively assumed to be discontinuous and given by 
$$q^\dag := 2\chi_{\Omega^q_1} + \chi^q_{\Omega^q_2} + 3\chi_{\Omega^q_3}$$
with 
$$\Omega^q_1:= (-1,-1/2)\times(-1,1), \quad \Omega^q_2:= (-1/2,1/2)\times(-1,1), \quad \Omega^q_2:= (1/2,1)\times(-1,1)$$
and
$$a^\dag := 3\chi_{\Omega^a_1} + 5\chi_{\Omega^a_2},$$
where
$$\Omega^a_1:= (-1,1)\times(-1,0), \quad \Omega^a_2:= (-1,1)\times(0,1).$$

The exact Neumann data on $\Gamma$ given by
\begin{align}\label{22-12-18ct1}
j^\dag := -\chi_{(0,1)\times\{-1\}} + \chi_{[-1,0]\times\{-1\}} -2\chi_{\{-1\}\times(-1,0]} + 3 \chi_{\{-1\}\times(0,1)}
\end{align}
and the exact Dirichlet data 
\begin{align}\label{5-6-20ct2}
g^\dag := \gamma_{|\Gamma} N_{j^\dag}(q^\dag,a^\dag).
\end{align}

The constants appearing in the sets $Q$ and $A$ are chosen as $\underline{q} := \underline{a} := 0.1$ and $\overline{q} := \overline{a} :=8$. The interval $(-1,1)$ is divided into $\tau$ equal segments and the domain $\Omega$ is then divided into $2\tau^2$ triangles with the diameter of each triangle $h = h_{\tau} = \sqrt{8}/{\tau}$. In the problem $\big(\mathcal{P}^h_{\delta,\rho}\big)$ the regularization parameter is taken by $\rho=\rho_\tau := 10^{-3}\sqrt{h}$ and the noisy observation data is assumed to be available in the form
\begin{align}\label{3-7-17ct1}
\left( j_{\delta_{\tau}}, g_{\delta_{\tau}} \right) = \left( j^\dag+ r\theta,~ g^\dag+ r\theta\right),
\end{align}
where $r$ is randomly generated in $(-1,1)$ and the positive parameter $\theta$ may depend on $\tau$ (cf.\ Example \ref{ex1}).

We utilize Algorithm \ref{31-8-18ct1} to reach the numerical solutions of the problem $\big(\mathcal{P}^h_{\delta,\rho}\big)$. The initial approximations are the constant functions defined by $q_0 = 1.5$ and $a_0 = 4$, the smoothing parameter $\epsilon^h = 10^{-3}\sqrt{h}$, positive constants $\kappa_1 = \kappa_2 = 10^{-3}\sqrt{h}$ and the maximum iterate $N=800$. The parameter $\beta$ appearing in \eqref{2-6-20ct3} is taken by $0.75$. We start the computational process with the coarsest level $\tau=4$ and then use the interpolation of the obtaining numerical solutions on the next finer mesh $\tau=8$ as initial approximations for the algorithm, and so on $\tau=16,32,64$.

With respect to the level $\tau$, we denote by $(q_\tau,a_\tau)$ the obtaining numerical solutions and then errors
\begin{alignat*}{2}
 E_{q,a} &= \big\|q_\tau -  q^\dag\big\|_{L^2(\Omega)} + \big\|a_\tau -  a^\dag\big\|_{L^2(\Omega)}, &\quad E_N &= \big\|N^{h_\tau}_{j_{\delta_\tau}} (q_\tau,a_\tau)  - N^{h_\tau}_{j^\dag}(q^\dag,a^\dag)\big\|_{L^2(\Omega)},\\	
 E_M &= \big\|M^{h_\tau}_{g_{\delta_\tau}} (q_\tau, a_\tau)  - M^{h_\tau}_{g^\dag}(q^\dag,a^\dag)\big\|_{L^2(\Omega)}, & E_D &= \Big\|D^{h_\tau}_{g^\dag_{\delta_\tau}} (q_\tau,a_\tau)  - D^{h_\tau}_{\widehat{g}^\dag}(q^\dag,a^\dag)\Big\|_{L^2(\Omega)},
\end{alignat*}
where 
$$g^\dag_{\delta_\tau} =
\begin{cases}
g_{\delta_\tau} & \mbox{on} \quad \Gamma,\\
\gamma_{|\partial\Omega\setminus\Gamma} N^{h_\tau}_{j^\dag}(q^\dag,a^\dag) & \mbox{on} \quad \partial\Omega\setminus \Gamma
\end{cases} \quad \mbox{and} \quad \widehat{g}^\dag := \gamma_{|\partial\Omega} N^{h_\tau}_{j^\dag}(q^\dag,a^\dag),$$
and $D^{h}_{g} (q,a)$ is the numerical solution of the problem
$
-\nabla \cdot \big(q \nabla u \big) + a u = f
$ 
in 
$\Omega$,
supplemented with the Dirichlet boundary condition $u = g$ on the boundary $\partial\Omega$.

\begin{example}\label{ex1}

To satisfy the condition \eqref{29-6-16ct1} we in this first implementation take $\theta$ in \eqref{3-7-17ct1} by
$$\theta = h_\tau \sqrt{10 \rho_\tau}.$$
and therefore the measurement noisy level is computed by
\begin{align}\label{5-6-20ct1}
\delta_\tau := \big\|j_{\delta_\tau} -j^\dag\big\|_{L^2(\Gamma)} + \big\|g_{\delta_\tau} -g^\dag\big\|_{L^2(\Gamma)}.
\end{align} 
The numerical result is summarized in Table \ref{b1}, where we present the different refinement levels $\tau$ and correspondingly noisy levels $\delta_\tau$ as well as the errors $E_{q,a},~ E_N,~ E_M$ and $E_D$. We observe that all errors and noisy levels get together smaller, as expected from our convergence result.

\begin{table}[H]
\begin{center}
\begin{tabular}{|l|l|l|l|l|l|}
\hline \multicolumn{6}{|c|}{ { Errors at refinement levels and correspondingly noisy levels} 
\vspace{0.05cm}
}\\

\hline
$\tau$ &\scriptsize $\delta_\tau$ &\scriptsize {$E_{q,a}$} &\scriptsize {$E_N$} &\scriptsize {$E_M$} &\scriptsize {$E_D$}\\
\hline
4   & 0.1912&  2.3025 & 0.5545& 0.3908 &0.1827\\
\hline
8  & 7.0192e-2&  0.7771& 0.1508& 0.1238 &6.9520e-2\\
\hline
16  & 2.6801e-2&  0.2712& 6.9911e-2& 4.9141e-2 &2.4563e-2\\
\hline
32  & 1.0663e-2&  0.1377& 3.5084e-2& 2.7613e-2 &1.6002e-2\\
\hline
64  & 4.3377e-3& 5.9782e-2& 1.6918e-2& 1.3122e-2 &8.0575e-3\\
\hline
\end{tabular}
\end{center}
\caption{Refinement levels $\tau$ and errors $E_{q,a}$, $E_N$, $E_M$, $E_D$ corresponding to the noise levels $\delta_\tau$.}
\label{b1}
\end{table}

Hereafter, all figures are presented with respect to the finest level $\tau = 64$. In Figure  \ref{h1} we from left to right show the graphs of the solution $q_\tau$ obtained from the computational process and the difference between the exact diffusion $q^\dag$ and the computed one $q_\tau$. The similar understanding for the reaction coefficient is presented in Figure \ref{h2}, meanwhile Figure \ref{h3} is utilized to perform the differences of the Neumann and Dirichlet boundary value problems, respectively.  

\begin{figure}[H]
\begin{center}
\includegraphics[scale=0.065]{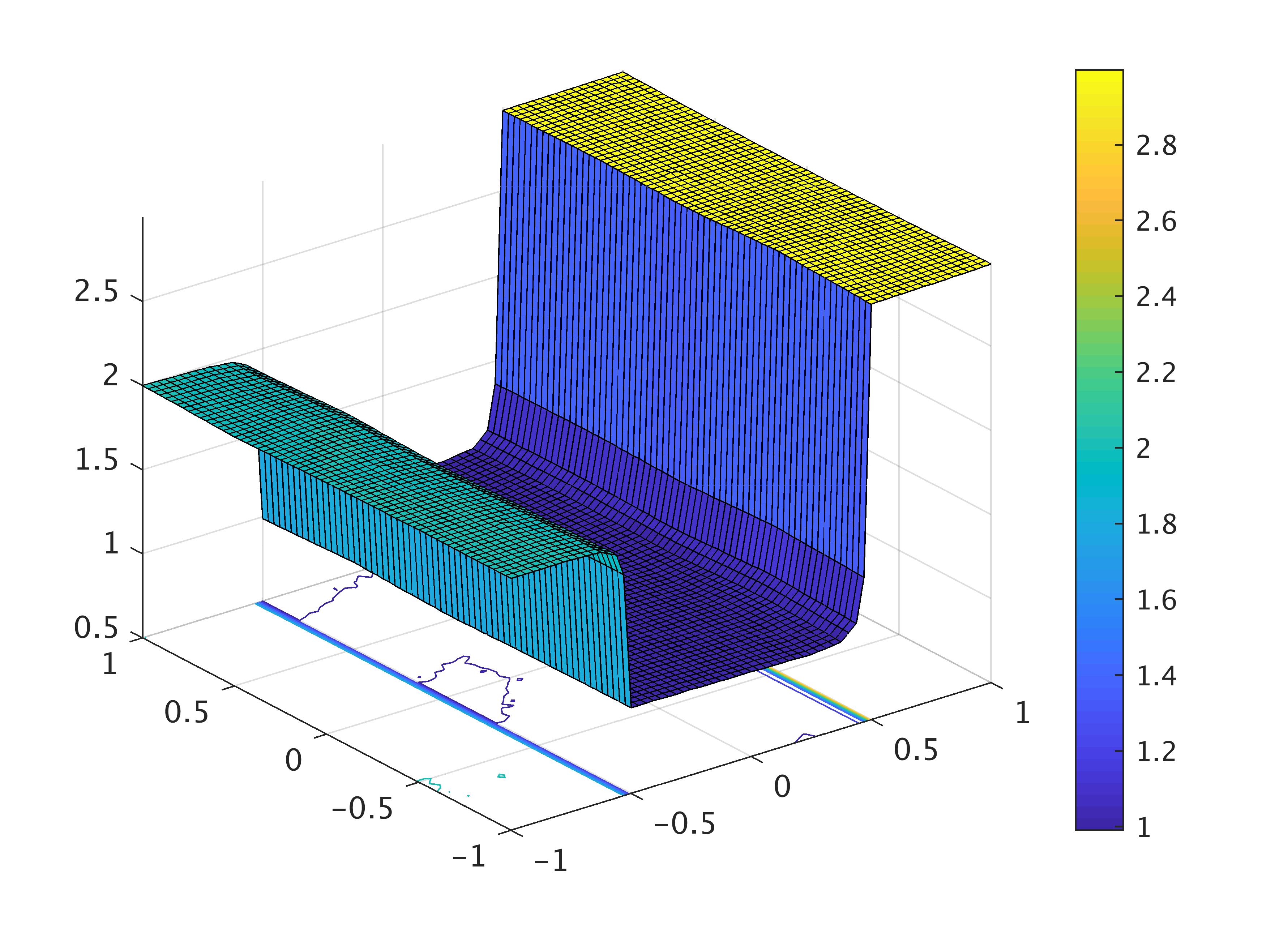}
\includegraphics[scale=0.065]{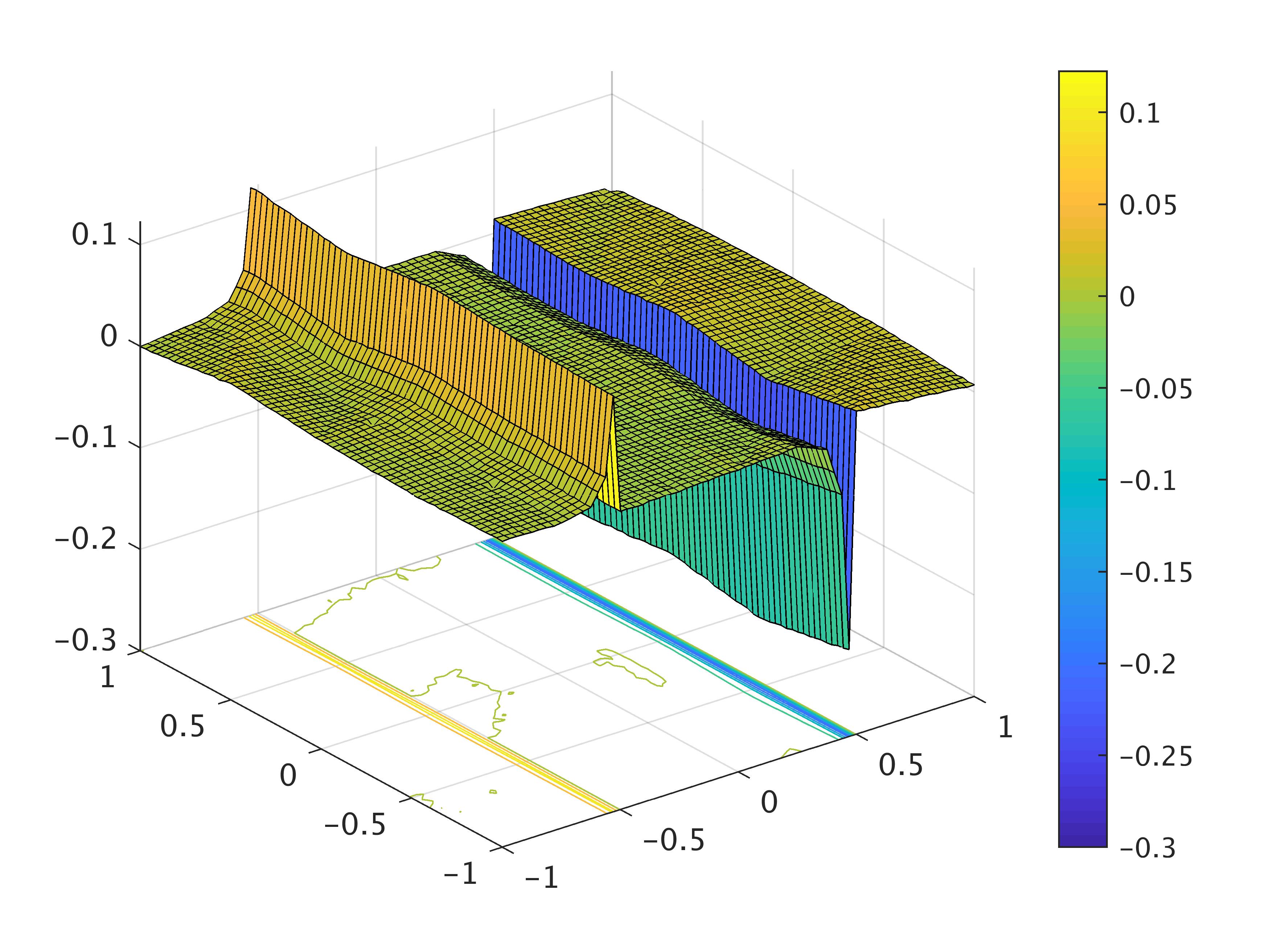} 
\end{center}
\caption{Computed diffusion $q_\tau$ (left) and the difference from the exact diffusion $q^\dag - q_\tau$ (right). }
\label{h1}
\end{figure}

\begin{figure}[H]
\begin{center}
\includegraphics[scale=0.065]{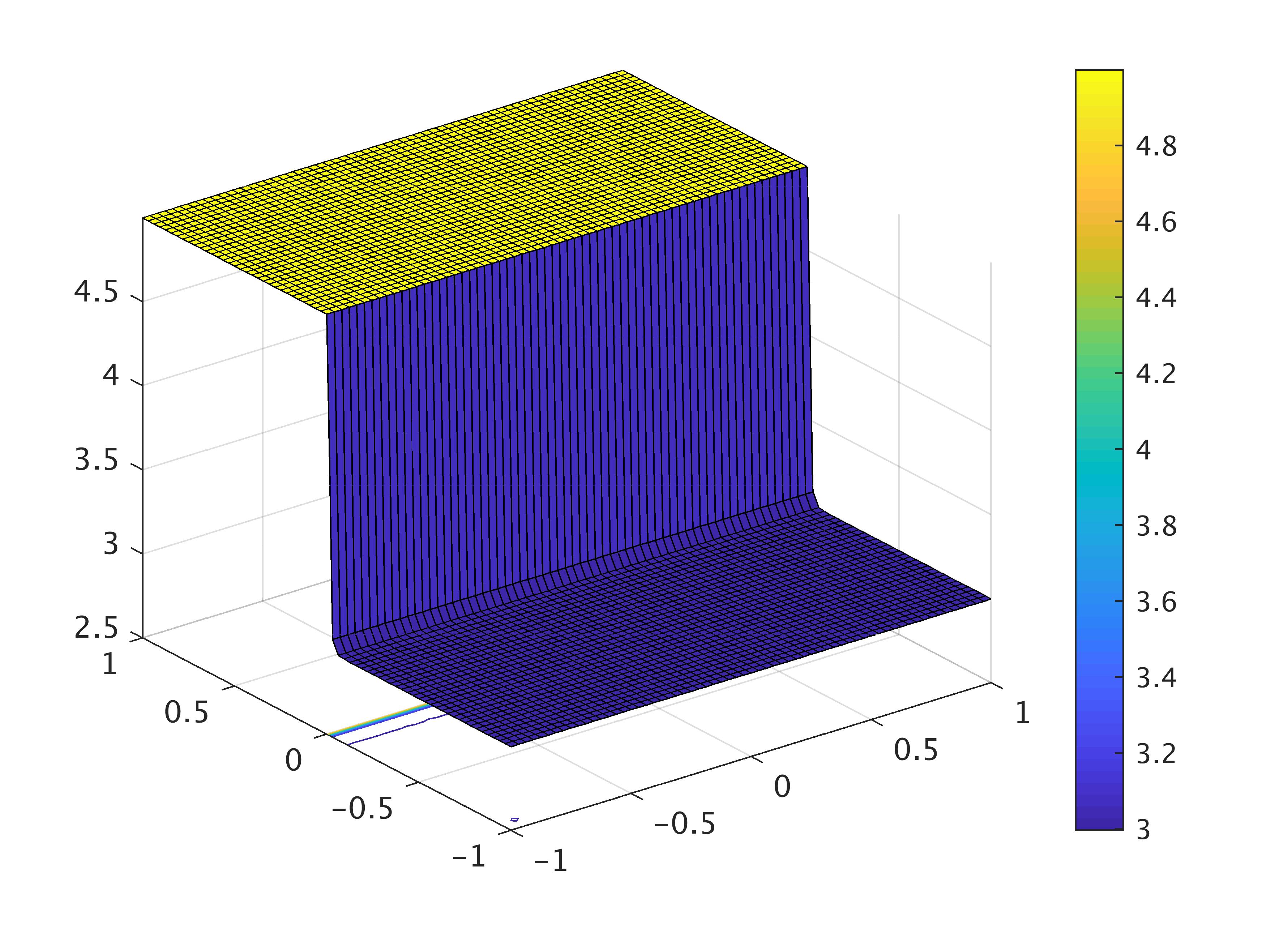}
\includegraphics[scale=0.065]{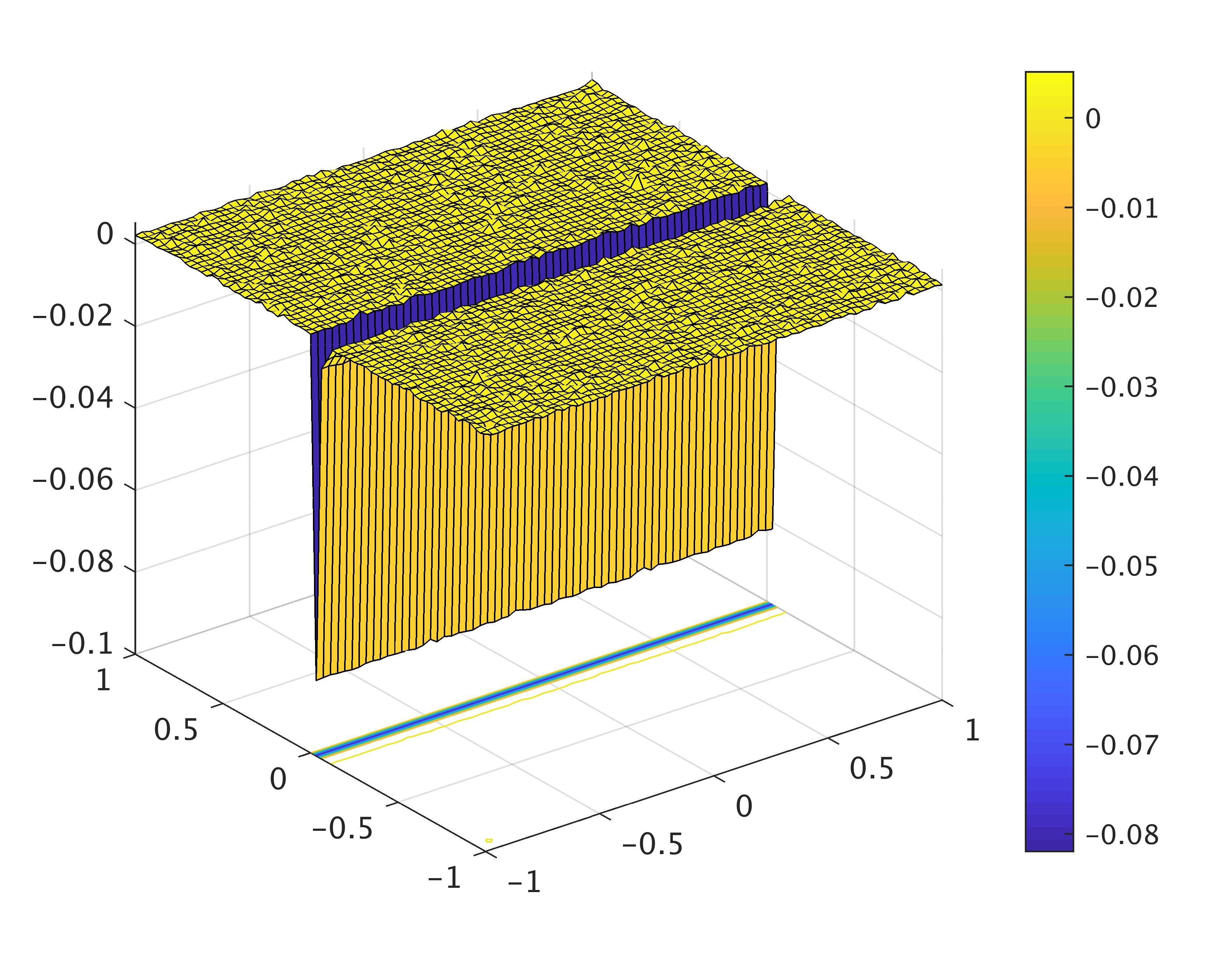} 
\end{center}
\caption{Computed reaction $a_\tau$ (left) and the difference from the exact reaction $a^\dag - a_\tau$ (right). }
\label{h2}
\end{figure}

\begin{figure}[H]
\begin{center}
\includegraphics[scale=0.065]{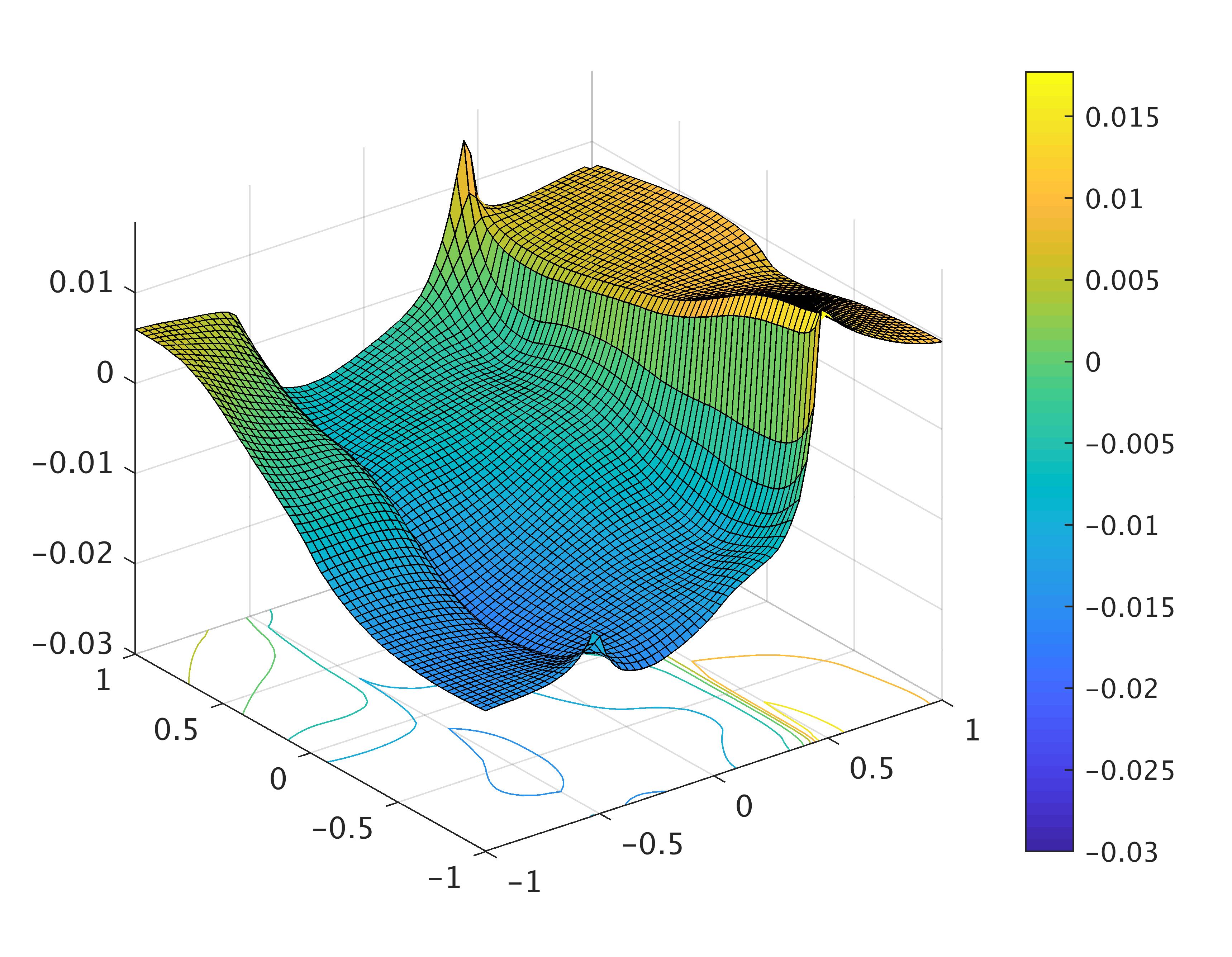}
\includegraphics[scale=0.065]{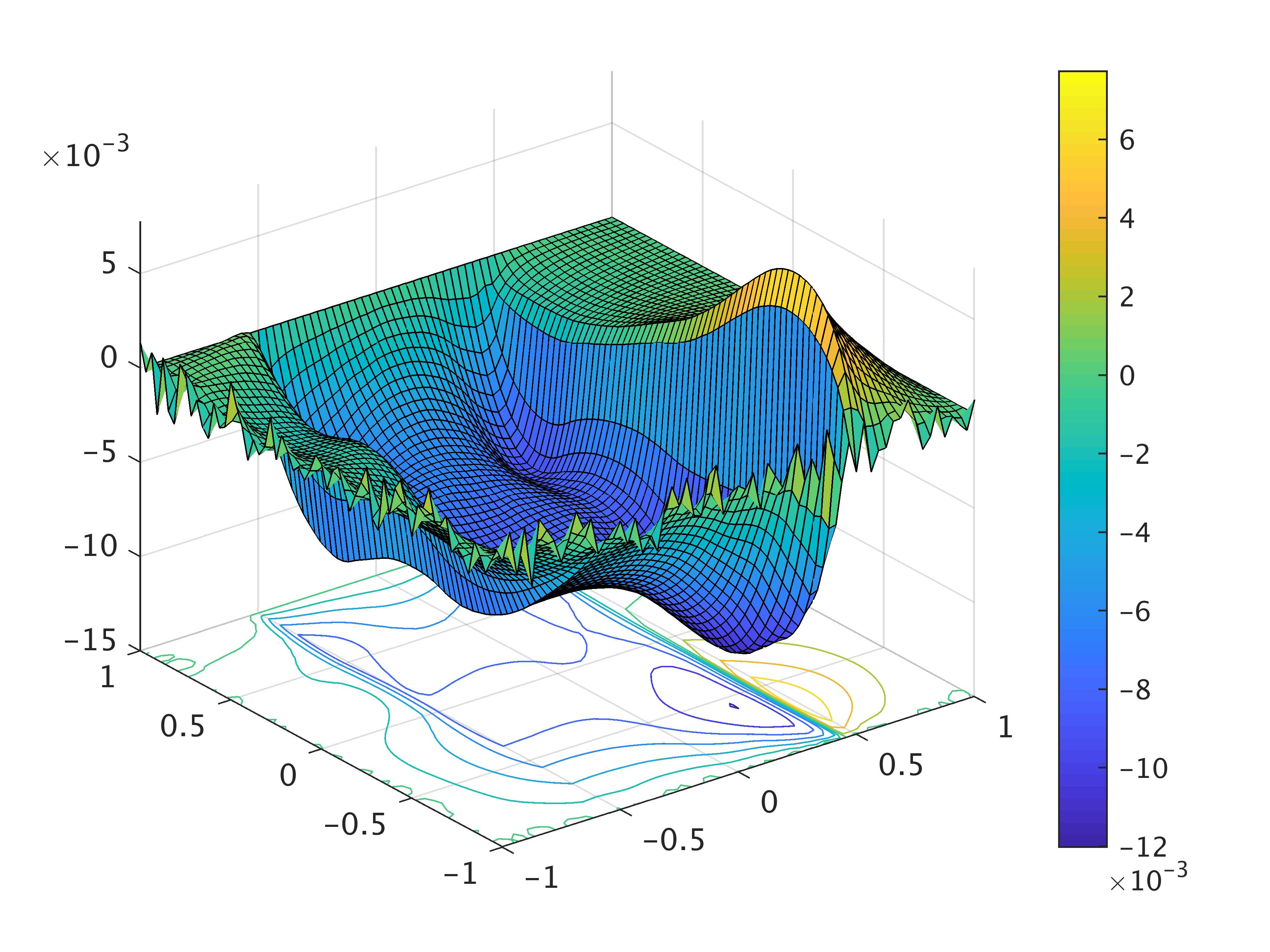} 
\end{center}
\caption{Differences $N^{h_\tau}_{j^\dag}(q^\dag,a^\dag) - N^{h_\tau}_{j_\delta}(q_\tau,a_\tau)$ (left) and $D^{h_\tau}_{\widehat{g}^\dag}(q^\dag,a^\dag) - D^{h_\tau}_{g^\dag_\delta}(q_\tau,a_\tau)$ (right). }
\label{h3}
\end{figure}

\end{example}


\begin{example}

We here take $\theta$ in \eqref{3-7-17ct1} to be  independent of the refinement levels $\tau$ by
$$\theta = 0.1, 0.05,~ \mbox{and}~ 0.01$$
and the noisy levels are then computed from the formula \eqref{5-6-20ct1}. In the Table \ref{rev-b1} we present the computational results for the finest refinement $\tau=64$.  As the previous implementation we observe a decrease of all errors as noisy levels get smaller.

\begin{table}[H]
\begin{center}
\begin{tabular}{|l|l|l|l|l|l|}
\hline \multicolumn{6}{|c|}{ {Numerical result for different values of $\theta$ at $\tau = 64$} }\\
\hline
$\theta$ &  $\delta$ &\scriptsize $E_{q,a}$ &\scriptsize $E_N$ &\scriptsize $E_M$ &\scriptsize $E_D$\\
\hline
0.01 & 3.0011e-2 & 9.4379e-2&  3.0362e-2 & 2.3549e-2 & 1.2186e-2 \\
\hline
0.05 & 0.1368 &   0.2115& 7.1326e-2& 5.0149e-2 & 2.6825e-2\\
\hline
0.1 & 0.2317 & 0.3921&  0.1485 & 0.1121 & 5.2141e-2\\
\hline
\end{tabular}
\end{center}
\caption{Numerical result for different values of $\theta$ independent of the refinement levels at the finest level $\tau=64$.}
\label{rev-b1}
\end{table}

We in the top two figures of Figure \ref{rev-h1} present the difference between the exact diffusion (respectively, reaction) and the computed diffusion (respectively, reaction) for $\theta=0.01$, while the similar illustrations for $\theta=0.05$ are performed in the bottom two figures.

\begin{figure}[H]
\begin{center}
\includegraphics[scale=0.065]{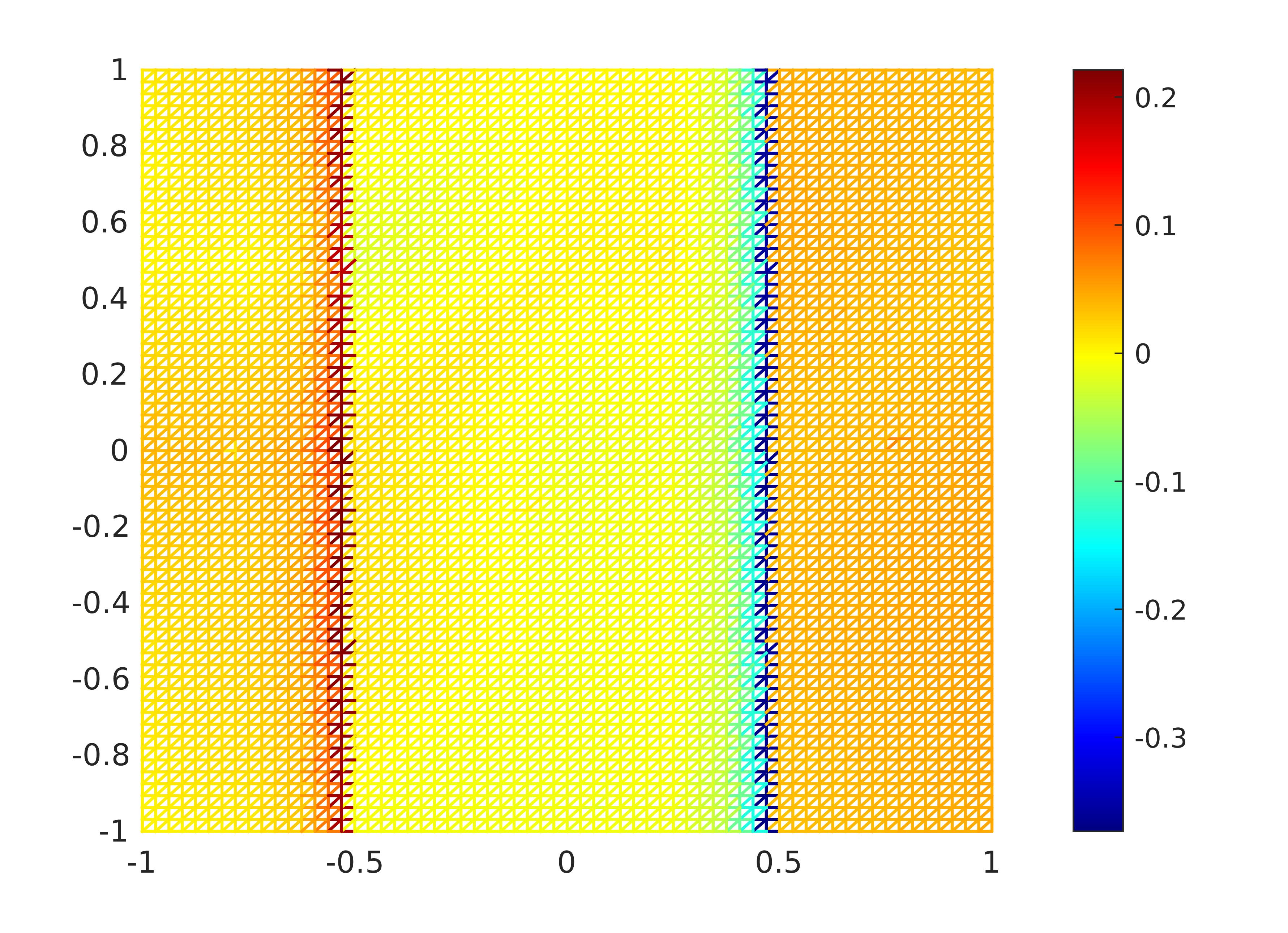}
\includegraphics[scale=0.065]{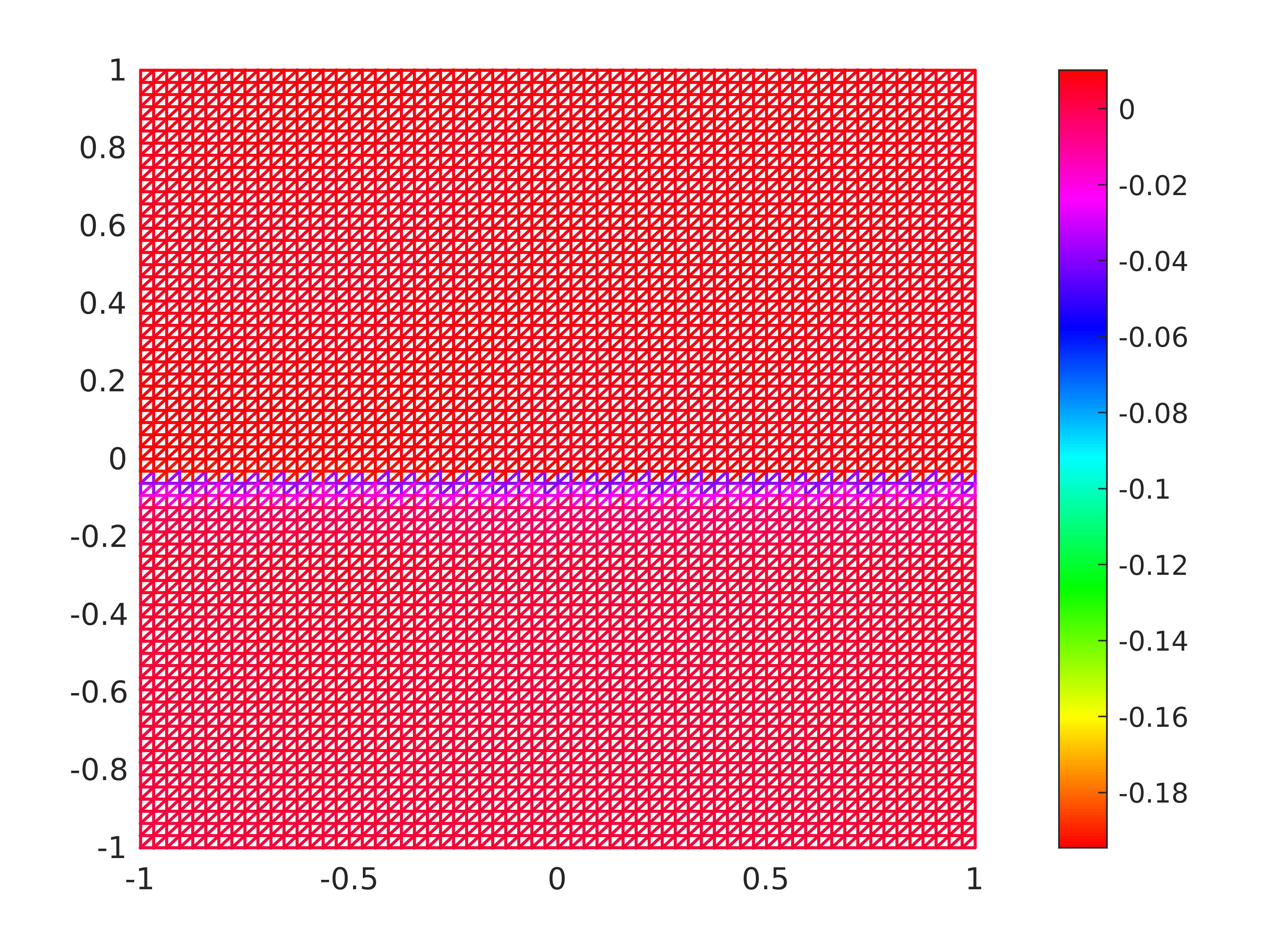}\\
\includegraphics[scale=0.065]{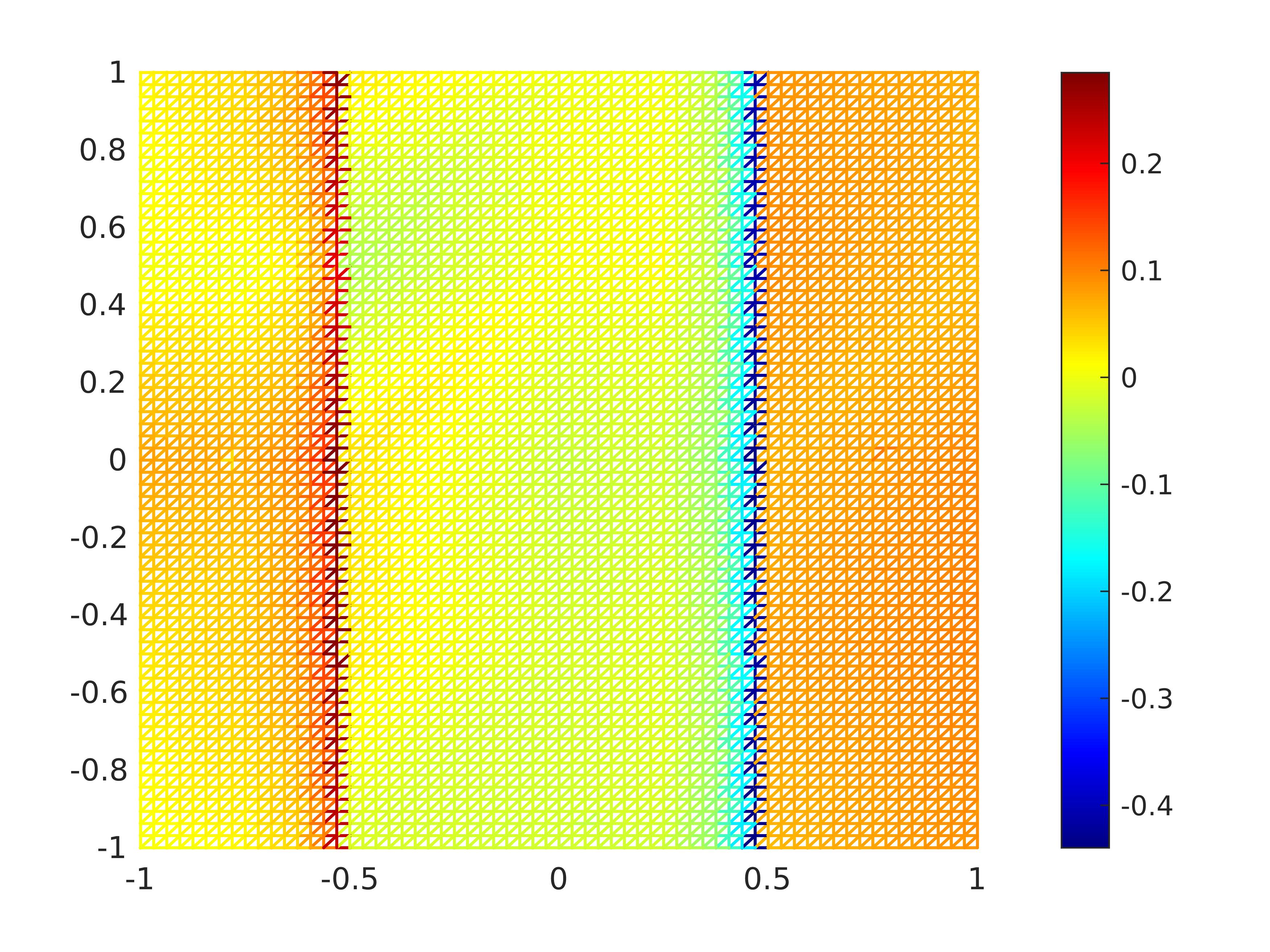}
\includegraphics[scale=0.065]{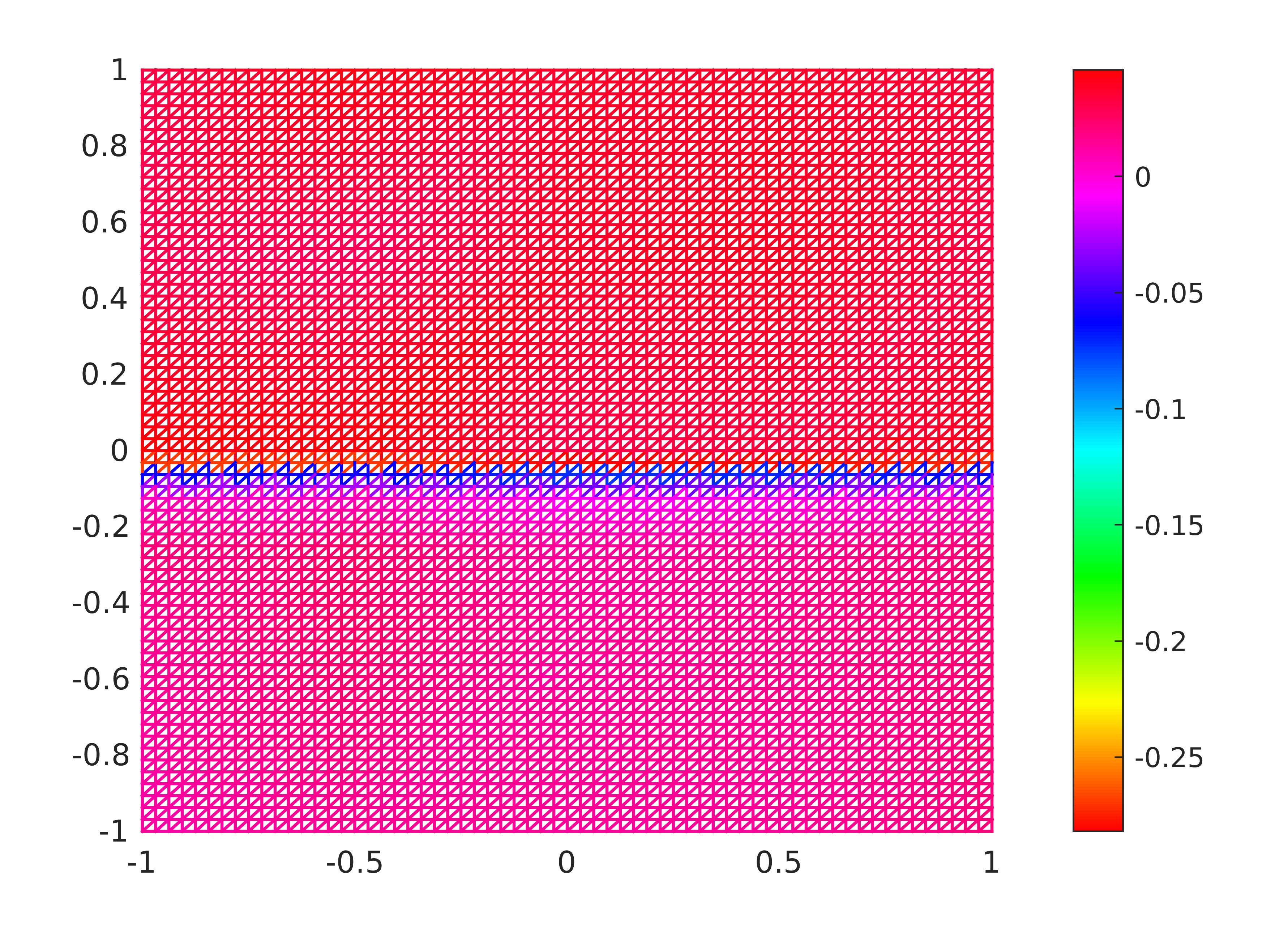}
\end{center}
\caption{Differences between the exact diffusion (respectively, reaction) and the computed diffusion (respectively, reaction): $\theta = 0.01$ (top two figures) and $\theta = 0.05$ (bottom two figures).}
\label{rev-h1}
\end{figure}

\end{example}


\begin{example}

Finally, we consider the case of multiple measurements. Assume that $I$ multiple measurements $\left(j_\delta^i,g_\delta^i \right)_{i=1,\ldots,I}$ on $\Gamma$  are available. With these datum at hand, we examine the minimization problem
$$\min_{(q,a)\in Q^h_{ad} \times A^h_{ad}} \widehat{\Upsilon}^h_{\delta,\rho} (q,a), \eqno \left(\widehat{\mathcal{P}}^h_{\delta,\rho}\right)$$
where
\begin{align*}
\widehat{\Upsilon}^h_{\delta,\rho} (q,a) &:=  \widehat{J}^h_\delta(q,a) + \rho R(q,a) \quad\mbox{and} \\
\widehat{J}^h_\delta(q,a) &:= \frac{1}{I}\sum_{i=1}^I \Bigg( \int_\Omega q\left|\nabla\big(N^h_{j^i_\delta}(q,a) - M^h_{g^i_\delta}(q,a)\big)\right|^2 dx + \int_\Omega a\big(N^h_{j^i_\delta}(q,a) - M^h_{g^i_\delta}(q,a)\big)^2 dx\\
&~\quad + \int_{\partial\Omega} \sigma\big(N^h_{j^i_\delta}(q,a) - M^h_{g^i_\delta}(q,a)\big)^2 ds \Bigg),
\end{align*}
which admits a minimizer $\big(\widehat{q}^h_{\delta,\rho}, \widehat{a}^h_{\delta,\rho}\big)$.

We now rewrite the exact boundary data $\left(j^\dag,g^\dag\right)$ in \eqref{22-12-18ct1} -- \eqref{5-6-20ct2} as
\begin{align*}
j^\dag_{(A,B,C,D)} := A \cdot \chi_{(0,1)\times\{-1\}} + B \cdot \chi_{[-1,0]\times\{-1\}} + C \cdot \chi_{\{-1\}\times(-1,0]} + D \cdot  \chi_{\{-1\}\times(0,1)}
\end{align*}
and $g^\dag_{(A,B,C,D)} :=  \gamma_{|\Gamma} N_{j^\dag_{(A,B,C,D)}} (q^\dag,a^\dag)$, which depend on the constants $(A,B,C,D)$. As in \eqref{3-7-17ct1}, the noisy observations are assumed to be given by
\begin{align}\label{15-6-20ct1}
\left( j^{(A,B,C,D)}_{\delta}, g^{(A,B,C,D)}_{\delta} \right) = \left( j^\dag_{(A,B,C,D)}+ r\theta, g^\dag_{(A,B,C,D)}+ r\theta\right).
\end{align}
In the case $(A,B,C,D)=(-1,1,-2,3)$ we obtain a single noisy measurement, i.e. $I=1$. We now fix $D=3$, and let $(A,B,C)$ take all permutations $\mathcal{S}_3$ of the set $\{-1,1,-2\}$, the equation \eqref{15-6-20ct1} then generates $I=6$ measurements. Likewise, if $(A,B,C,D)$ takes all permutations $\mathcal{S}_4$ of $\{-1,1,-2,3\}$ we have $I=16$ measurements.

The noisy level is given by
$$
\delta = \begin{cases}
 \big\|j^{(-1,1,-2,3)}_{\delta} -j^\dag_{(-1,1,-2,3)}\big\|_{L^2(\Gamma)} + \big\|g^{(-1,1,-2,3)}_{\delta} -g^\dag_{(-1,1,-2,3)}\big\|_{L^2(\Gamma)} \mbox{~if~}  (A,B,C,D)=(-1,1,-2,3),\\
\dfrac{1}{6}\sum_{(A,B,C)\in\mathcal{S}_3} \big\|j^{(A,B,C,3)}_{\delta} -j^\dag_{(A,B,C,3)}\big\|_{L^2(\Gamma)} + \big\|g^{(A,B,C,3)}_{\delta} -g^\dag_{(A,B,C,3)}\big\|_{L^2(\Gamma)} \mbox{~if~}  D=3,\\
\dfrac{1}{16}\sum_{(A,B,C,D)\in\mathcal{S}_4} \big\|j^{(A,B,C,D)}_{\delta} -j^\dag_{(A,B,C,D)}\big\|_{L^2(\Gamma)} + \big\|g^{(A,B,C,D)}_{\delta} -g^\dag_{(A,B,C,D)}\big\|_{L^2(\Gamma)}.
\end{cases}$$
A computation with $\theta=0.1$ and $\tau = 64$ shows $\delta = 0.2317$.
The correspondingly numerical result for the multiple measurement case is presented in the Table \ref{b2}, where its first line is copied from the last one of Table \ref{rev-b1}. We observe that the use of multiple measurements improves the obtaining numerical solutions in case of the large noise level, as can be seen all errors decrease oppositely with the increase of the number of measurements.

\begin{table}[H]
\begin{center}
\begin{tabular}{|c|l|l|l|l|}
\hline \multicolumn{5}{|c|}{ {Errors for multiple measurements} }\\
\hline
\mbox{Number of measurements} $I$  &\scriptsize $E_{q,a}$ &\scriptsize $E_N$ &\scriptsize $E_M$ &\scriptsize $E_D$\\
\hline
1  & 0.3921&  0.1485 & 0.1121 & 5.2141e-2 \\
\hline
6 &   0.2849& 8.4842e-2& 6.0352e-2 & 3.8542e-2\\
\hline
16 &  0.1609& 4.7708e-2& 3.3418e-2 & 1.9260e-2\\
\hline
\end{tabular}
\end{center}
\caption{Errors for multiple measurements $I=1,6,16$ with $\theta = 0.1$ and $\tau = 64$, i.e. $\delta = 0.2317$.}
\label{b2}
\end{table}

Finally, in Figure \ref{h4} -- Figure \ref{h6} we perform the graphs of the computation, which include the  differences between the exact coefficient and the computational one, i.e. $q^\dag - q_\tau$ (left) and $a^\dag - a_\tau$ (right). We would like to note that the computational errors occur much more at areas where the identified coefficients are discontinuous than others.

\begin{figure}[H]
\begin{center}
\includegraphics[scale=0.065]{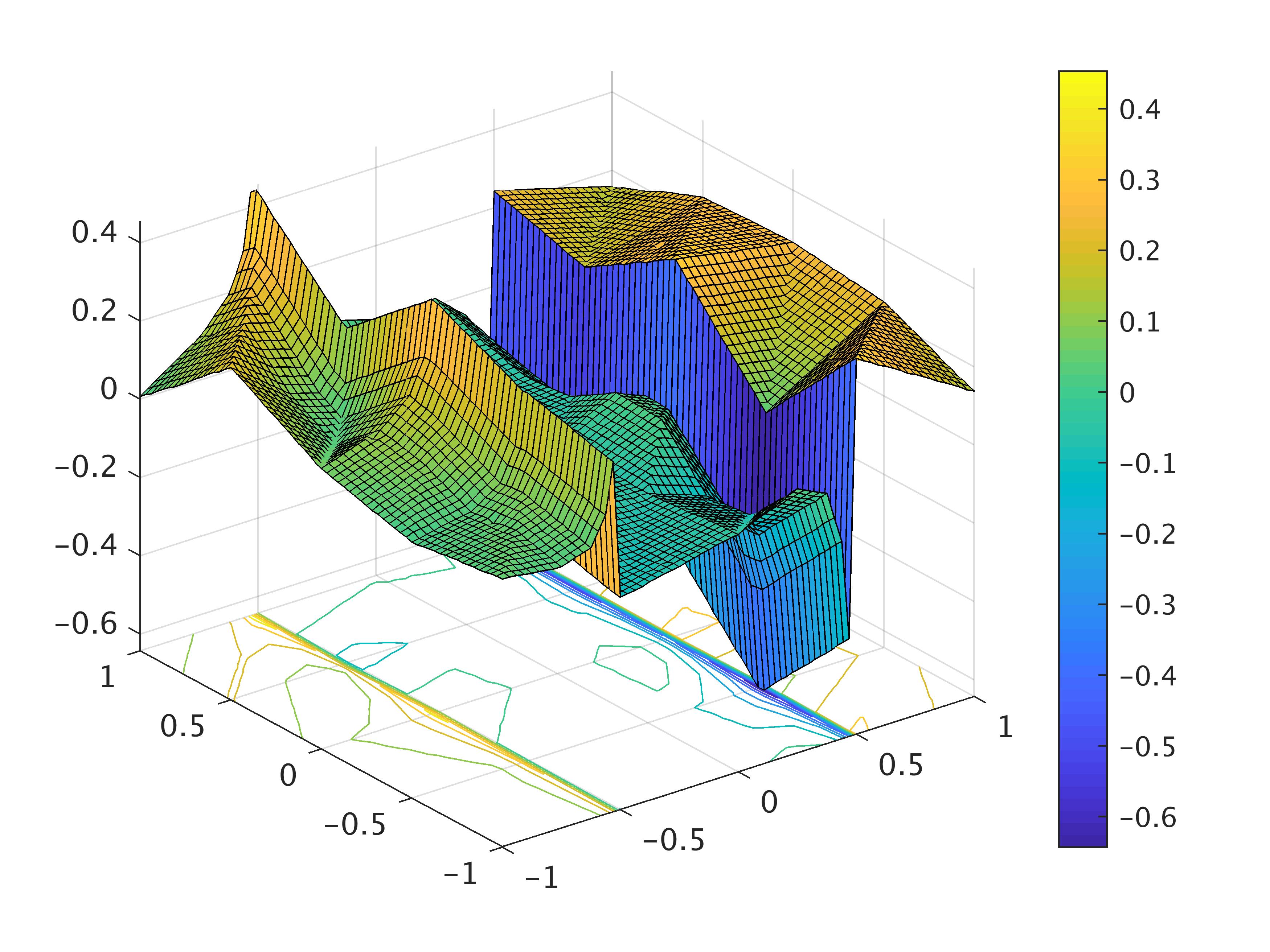}
\includegraphics[scale=0.065]{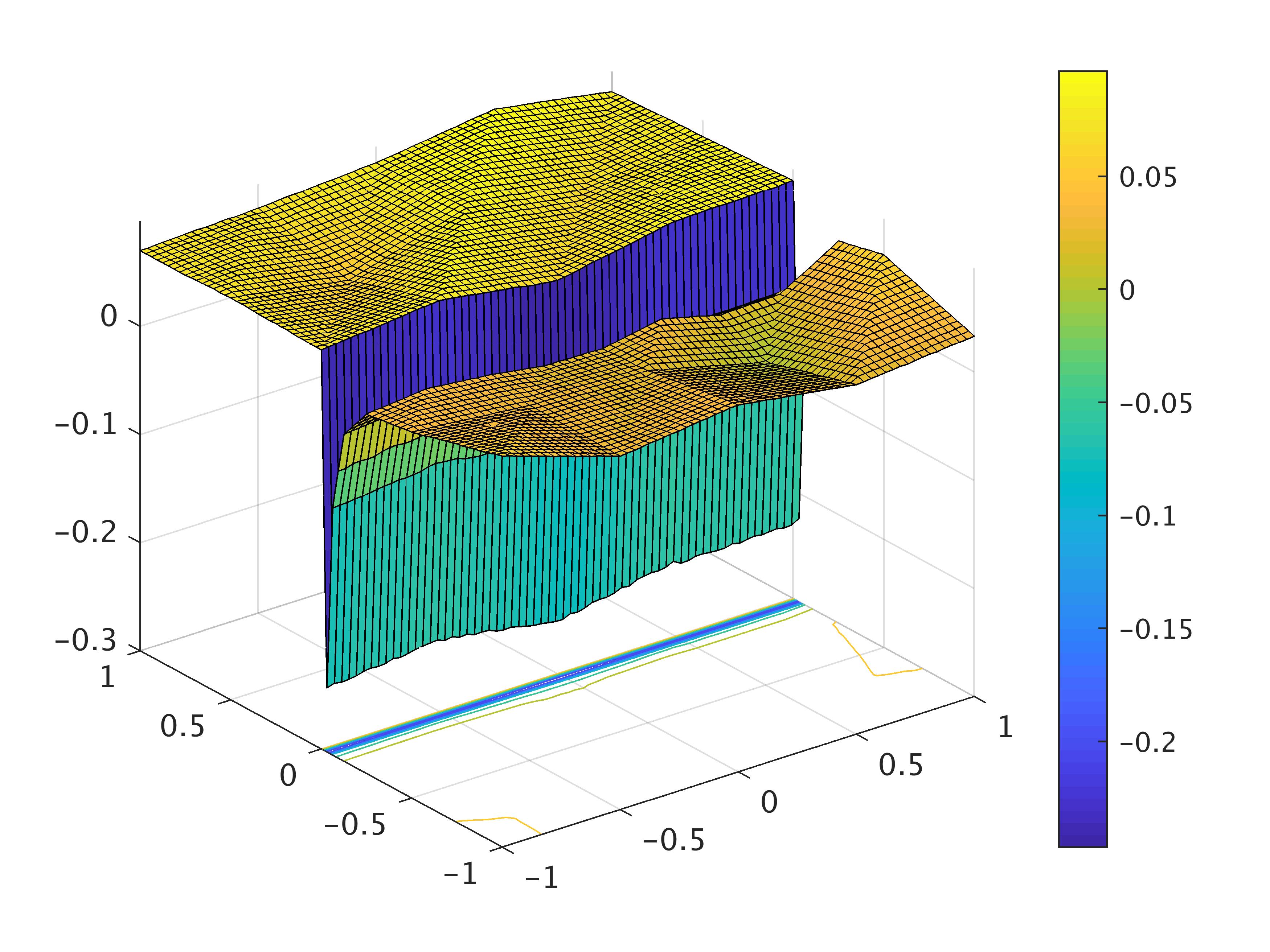} 
\end{center}
\caption{$I=1$: differences $q^\dag - q_\tau$ (left) and $a^\dag - a_\tau$ (right).}
\label{h4}
\end{figure}

\begin{figure}[H]
\begin{center}
\includegraphics[scale=0.065]{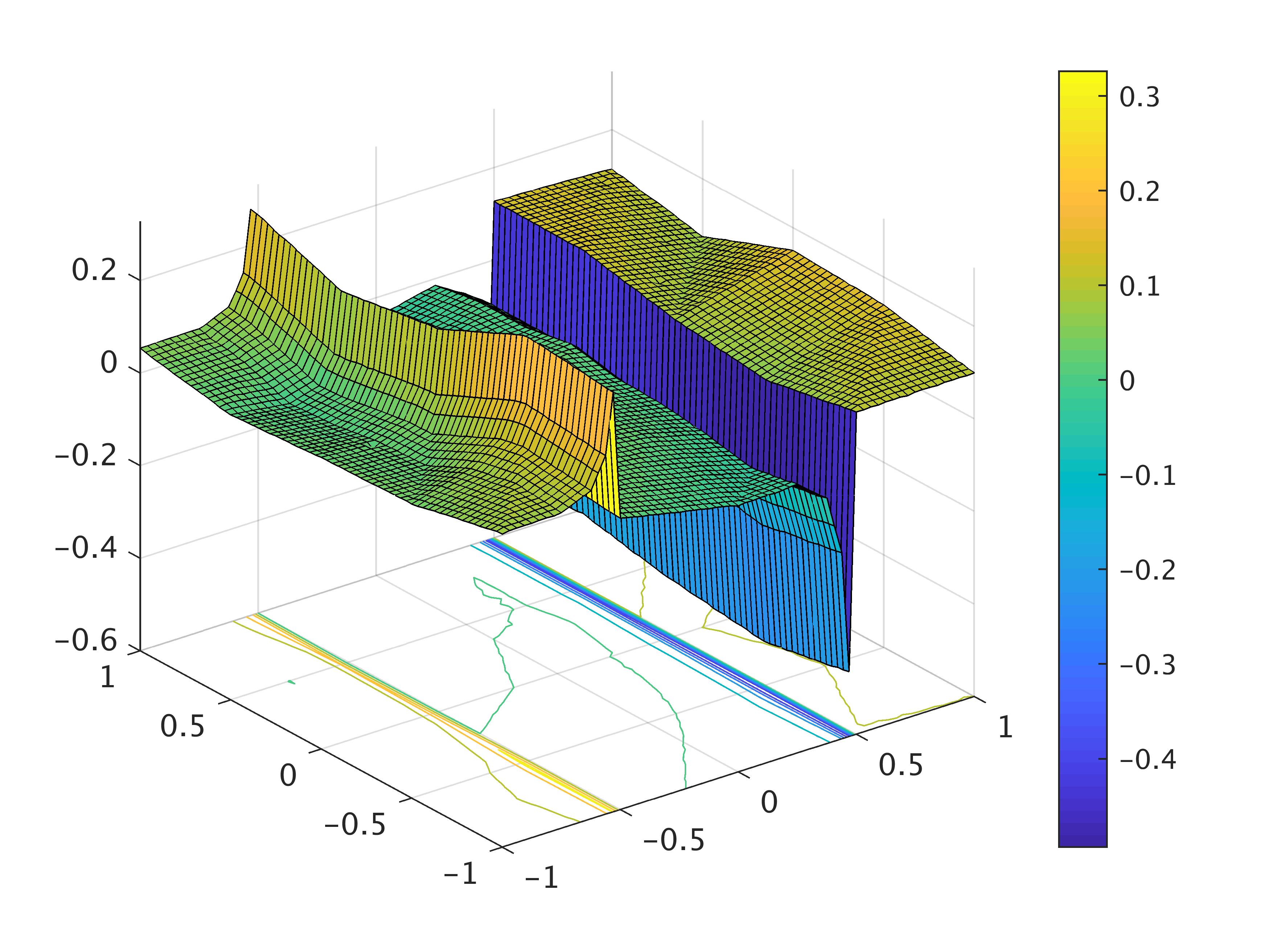}
\includegraphics[scale=0.065]{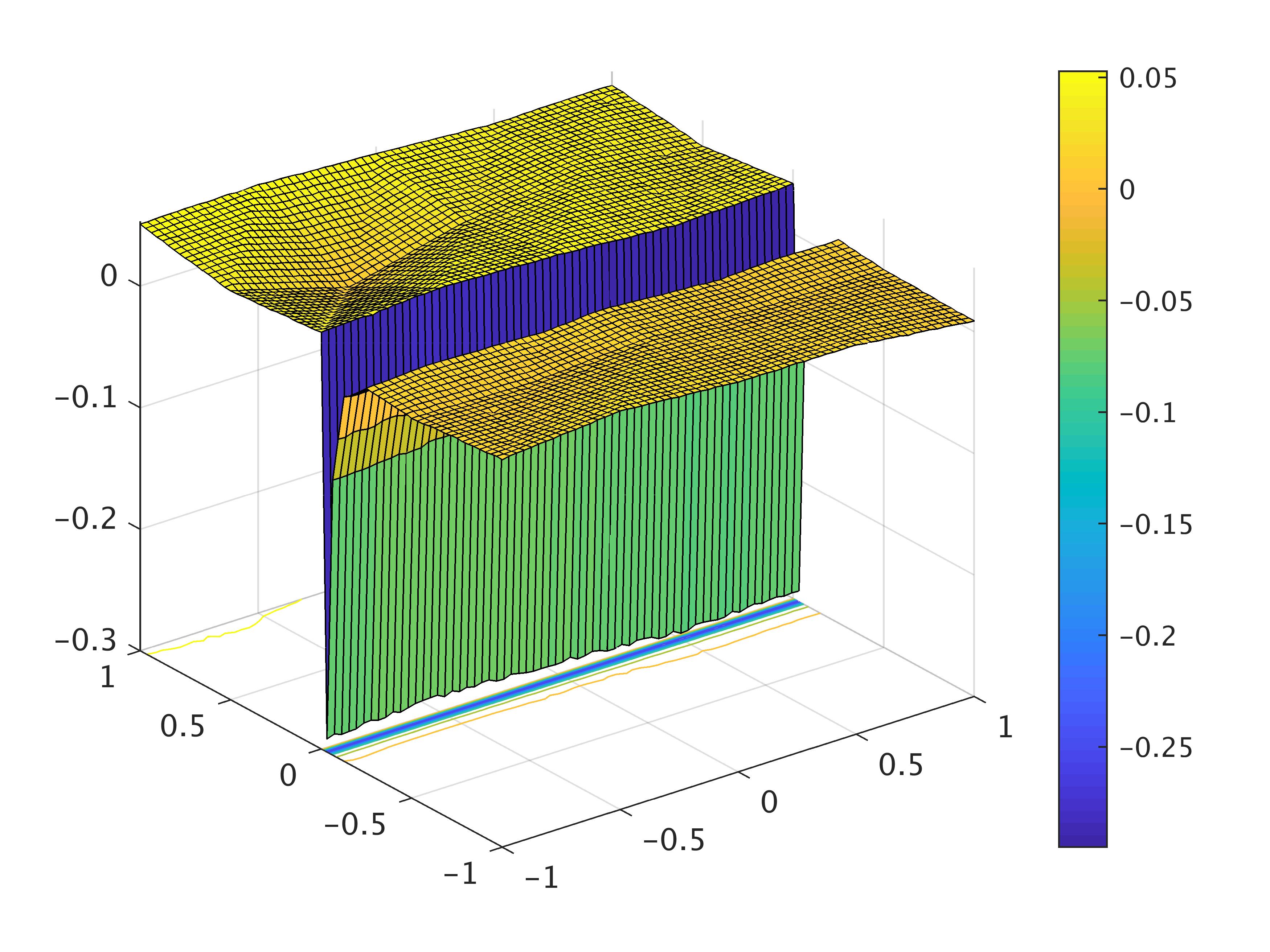} 
\end{center}
\caption{$I=6$: differences $q^\dag - q_\tau$ (left) and $a^\dag - a_\tau$ (right).}
\label{h5}
\end{figure}

\begin{figure}[H]
\begin{center}
\includegraphics[scale=0.065]{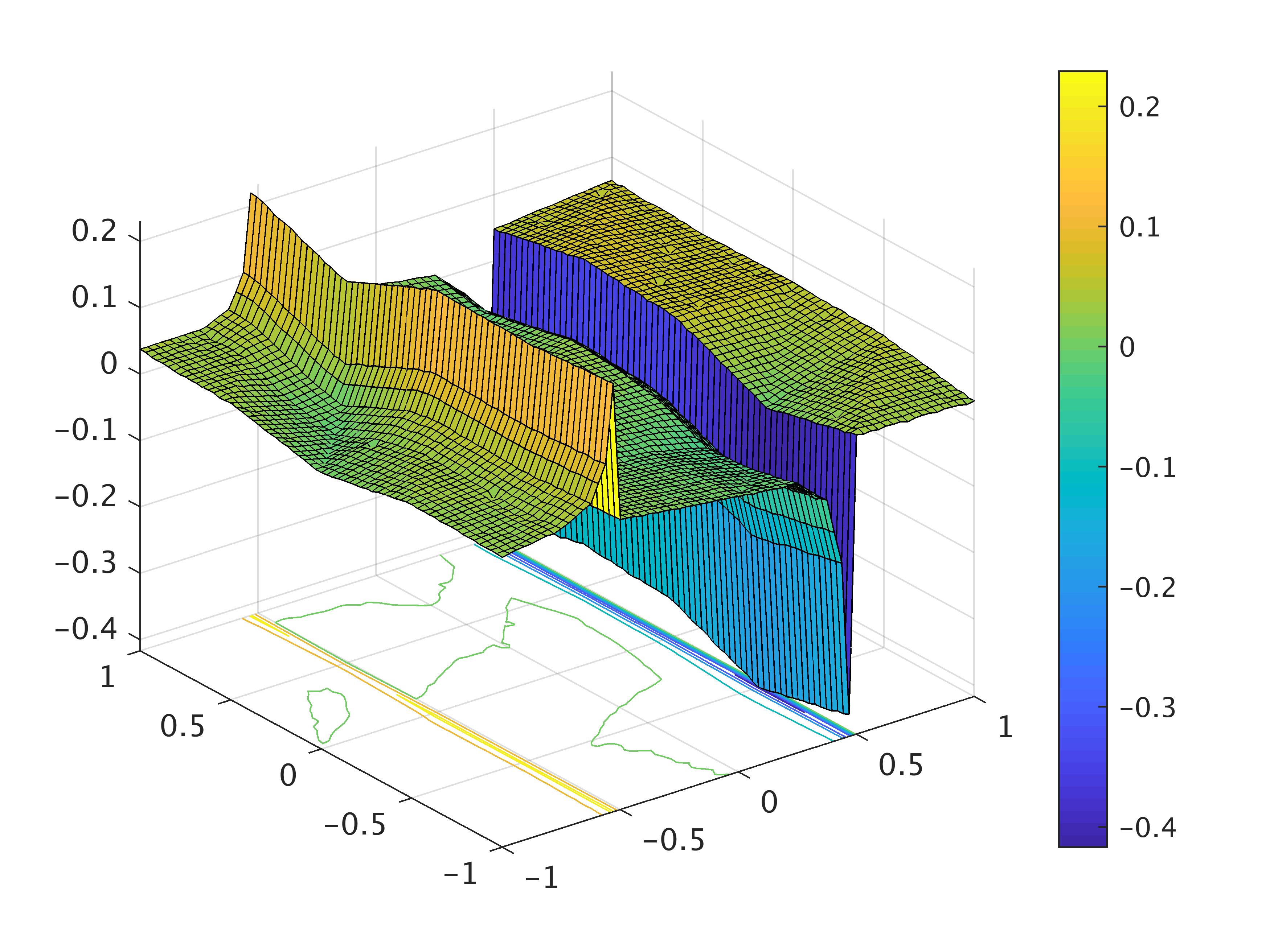}
\includegraphics[scale=0.065]{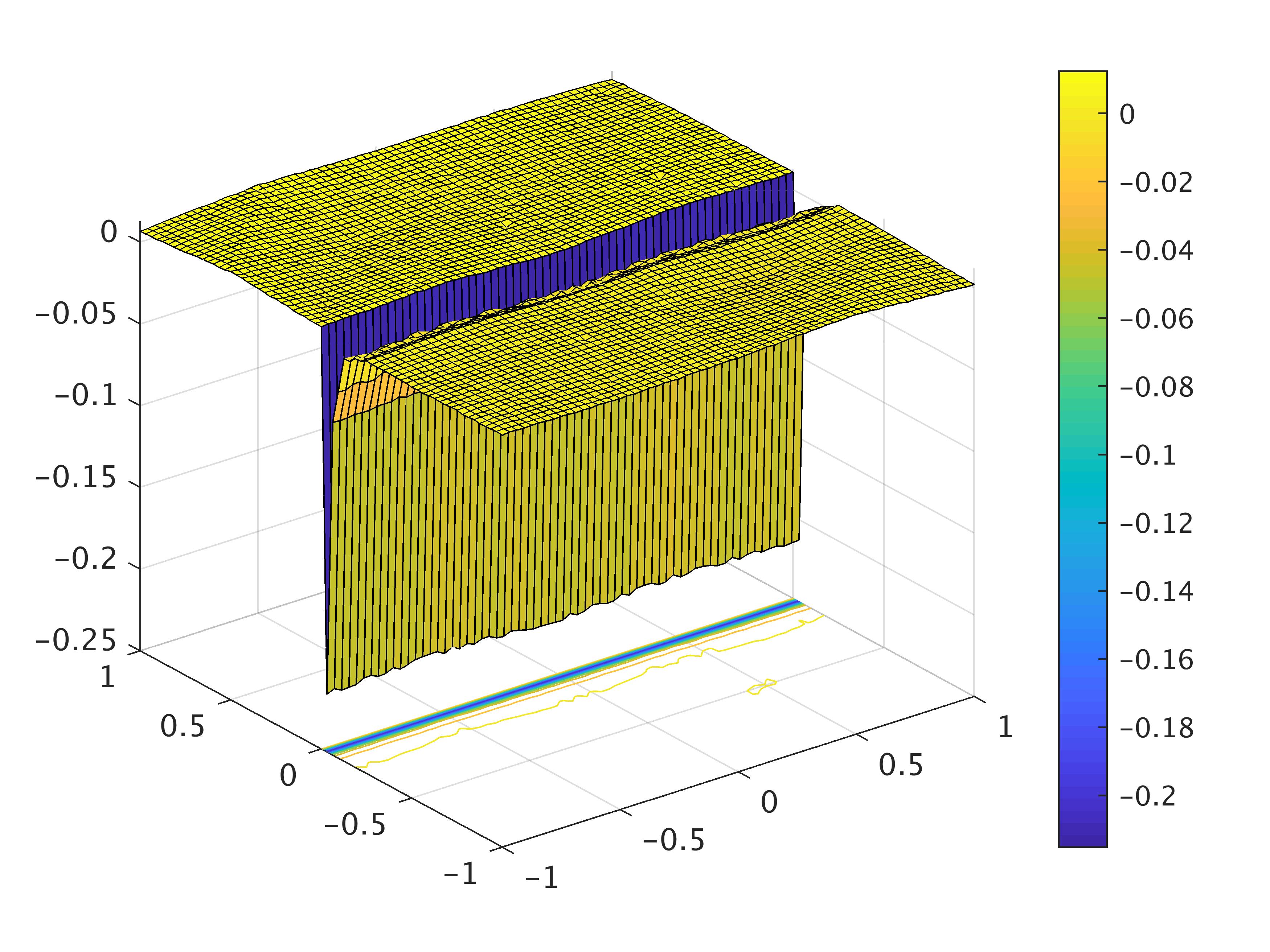} 
\end{center}
\caption{$I=16$: differences $q^\dag - q_\tau$ (left) and $a^\dag - a_\tau$ (right).}
\label{h6}
\end{figure}

\end{example}


{\bf Acknowledgements:} The author would like to thank the referees and the editor for their valuable comments and suggestions which helped to improve our paper.


\end{document}